\documentclass[11pt]{article}
\usepackage[francais]{babel}
\usepackage{amsmath, amssymb, amsfonts, amsthm, amscd, vmargin}
\usepackage[latin1]{inputenc}
\usepackage[all]{xy}
\usepackage{graphicx}

\setmarginsrb{2.4cm}{4.5cm}{2.4cm}{3.3cm}{0cm}{0cm}{0cm}{1.5cm}

\theoremstyle{plain}
\newtheorem{teo}{Th\'eor\`eme}[section]
\newtheorem{lem}[teo]{Lemme}
\newtheorem{prop}[teo]{Proposition}
\newtheorem{cor}[teo]{Corollaire}

\theoremstyle{definition}
\newtheorem{defi}[teo]{D\'efinition}
\newtheorem{ex}[teo]{Exemple}
\newtheorem{exs}[teo]{Exemples}
\newtheorem{rem}[teo]{Remarque}
\newtheorem{rems}[teo]{Remarques}
\numberwithin{equation}{teo}

\renewcommand{\dim}{\operatorname{dim}}

\newcommand{\Cbb}{{\mathbb C}}
\newcommand{\Qbb}{{\mathbb Q}}
\newcommand{\Zbb}{{\mathbb Z}}
\newcommand{\Rbb}{{\mathbb R}}
\newcommand{\Pbb}{{\mathbb P}}

\newcommand{\lra}{\longrightarrow}
\newcommand{\lmt}{\longmapsto}

\begin{document}

\title{Vari\'et\'es horosph\'eriques de Fano}

\author{Boris Pasquier}
\maketitle

\section*{Introduction}

Une vari\'et\'e complexe projective $X$ est dite de Fano si elle est normale et si son diviseur anticanonique $-K_X$ est de Cartier et ample. On sait qu'il existe seulement un nombre fini de familles de vari\'et\'es lisses, de Fano et de dimension donn\'ee. Cependant ces familles sont seulement connues jusqu'en dimension~$3$.

Les vari\'et\'es toriques donnent beaucoup d'exemples de vari\'et\'es de Fano. Plus pr\'ecis\'ement, V.~Batyrev a classifi\'e les vari\'et\'es toriques de Fano de dimension~$n$ en termes des polytopes r\'eflexifs de dimension~$n$ \cite{Ba94}: ce sont les polytopes convexes de $\Rbb^n$ \`a sommets dans $\Zbb^n$ contenant~$0$ dans leur int\'erieur et tels que leur polytope dual v\'erifie les m\^emes hypoth\`eses.

De plus, certaines propri\'et\'es ou certains invariants g\'eom\'etriques des vari\'e\-t\'es toriques de Fano, comme la lissit\'e, le nombre de Picard ou le degr\'e, se lisent facilement sur le polytope r\'eflexif associ\'e. Cela a permis \`a O.~Debarre de majorer le degr\'e $(-K_X)^d$ des vari\'et\'es toriques lisses de Fano en fonction de la dimension~$d$ et du nombre de Picard \cite{De03}.  D'autre part, C.~Casagrande a r\'ecemment donn\'e une majoration optimale du nombre de Picard des vari\'et\'es toriques $\Qbb$-factorielles et de Fano en fonction de la dimension \cite{Ca06}.

Cet article a pour but de g\'en\'eraliser tous ces r\'esultats aux vari\'et\'es horosph\'eriques. Soit $G$ un groupe alg\'ebrique r\'eductif connexe. Un $G$-espace homog\`ene  est dit horosph\'erique de rang $n$ si c'est un fibr\'e en tores $(\Cbb^*)^n$ sur une vari\'et\'e de drapeaux. Voici quelques exemples d'espaces homog\`enes horosph\'eriques $G/H$ :\\
\begin{center}
\begin{tabular}{|c|c|c|c|c|}
\hline
 & $G$ & $H$ &  rang & dimension\\
\hline
\hline
$1$ & $(\Cbb^*)^n$ & $\{1\}$ & $n$ & $n$\\
\hline
$2$ & $G$ &  un sous-groupe parabolique $P$ & $0$ & $\dim G-\dim P$\\
\hline
$3$ & $SL_2$ & $U=\{\left(\begin{array}{cc}1&*\\0&1\end{array}\right)\}$ & $1$ & $2$\\
\hline
$4$ & $SL_2\times\Cbb^*$ & $U=\{\left(\begin{array}{cc}1&*\\0&1\end{array}\right)\}\times\{1\}$ & $2$ & $3$ \\
\hline
$5$ & $SL_2\times SL_2$ & $U=\{\left(\begin{array}{cc}1&*\\0&1\end{array}\right)\}\times\{\left(\begin{array}{cc}1&*\\0&1\end{array}\right)\}$ & $2$ & $4$  \\
\hline
$6$ & $SL_3$ & $U=\{\left(\begin{array}{ccc}1&*&*\\0&1&*\\0&0&1\end{array}\right)\}$ & $2$ & $5$ \\
\hline
\end{tabular}
\end{center}

Une vari\'et\'e horosph\'erique est un plongement d'un espace homog\`ene horosph\'erique $G/H$, c'est-\`a-dire une $G$-vari\'et\'e normale contenant une orbite ouverte isomorphe \`a $G/H$; son rang est celui de $G/H$.  Parmi les vari\'et\'es horosph\'eriques, on compte les vari\'et\'es toriques (lorsque $G/H$ est un tore: exemple $1$) et les vari\'et\'es de drapeaux (exemple $2$). Ces derni\`eres sont lisses et de Fano.

Les vari\'et\'es horosph\'eriques font partie de la famille des vari\'et\'es sph\'eriques. Les plongements d'un espace homog\`ene sph\'erique $G/H$ fix\'e ont \'et\'e classifi\'es en termes d'\'eventails colori\'es par D.~Luna et T.~Vust \cite{LV83}.
Lorsque $G/H$ est horosph\'erique de rang $n$, on montre que les plongements de Fano de $G/H$ sont classifi\'es en termes de certains polytopes rationnels, dits $G/H$-r\'eflexifs (voir la d\'efinition \ref{reflexif}). Ces polytopes sont de dimension~$n$ (tout comme les \'eventails colori\'es). Il est important de remarquer que la dimension de $G/H$ est plus grande que $n$, avec \'egalit\'e si et seulement si $G/H$ est un tore; dans ce  dernier cas, les polytopes $G/H$-r\'eflexifs sont les polytopes r\'eflexifs d\'efinis par V.~Batyrev. A rang \'egal, les polytopes $G/H$-r\'eflexifs peuvent \^etre beaucoup plus nombreux que les polytopes r\'eflexifs. Lorsque $G$ et $H$ sont comme dans l'exemple $6$, il y a, \`a automorphisme pr\`es, $398$ polytopes $G/H$-r\'eflexifs \cite[ch.6]{Pa06}. En comparaison, on compte seulement $16$ polytopes r\'eflexifs de dimension~$2$.

V.~Alexeev et M.~Brion ont montr\'e que l'ensemble des classes d'isomorphisme des vari\'et\'es sph\'eriques de Fano de dimension fix\'ee est fini \cite{AB04}. On verra que la classification pr\'ec\'edente permet d'avoir une version effective de ce r\'esultat pour les vari\'et\'es horosph\'eriques de Fano dont l'orbite ouverte est fix\'ee.\\

Dans la partie~1, on pr\'esente la classification de Luna et Vust dans le cas d'un espace homog\`ene horosph\'erique $G/H$.

Un crit\`ere de lissit\'e est donn\'e pour les vari\'et\'es horosph\'eriques dans la partie~2. Ce crit\`ere, qui g\'en\'eralise le r\'esultat de F.~Pauer \cite{Pa83}, a aussi \'et\'e obtenu r\'ecemment par  D.~Timashev \cite[th. 28.3]{Ti06}. On montre aussi que, comme dans le cas torique, toute sous-vari\'et\'e irr\'eductible et stable par $G$ d'une vari\'et\'e horosph\'erique lisse, est aussi lisse.

Dans la partie~3, on  classifie les plongements de Fano de $G/H$ en termes de polytopes $G/H$-r\'eflexifs, et on donne une borne explicite du nombre de classes d'isomorphisme de plongements de Fano de $G/H$.

Gr\^ace \`a cette classification, on d\'emontre les r\'esultats suivants, dans la partie~4.
\begin{teo}\label{majorationdegre}
Soit $X$ une vari\'et\'e horosph\'erique de Fano, localement factorielle, de dimension $d$, de rang $n$ et de nombre de Picard $\rho$.\\
Si  $\rho>1$ alors $$(-K_X)^d\leq d!\,d^{d\rho+n}.$$
Si $\rho=1$, on a $$(-K_X)^d\leq d!\,(d+1)^{d+n}.$$
\end{teo}
Remarquons qu'une vari\'et\'e lisse est toujours localement factorielle. La r\'eciproque est vraie pour les vari\'et\'es toriques mais elle est fausse pour les vari\'et\'es horosph\'eriques.
\begin{teo}\label{nombrepicard}
Soit $X$ une vari\'et\'e horosph\'erique de Fano, $\Qbb$-factorielle, de dimension $d$, de rang $n$ et de nombre de Picard $\rho$. On a $$\rho\leq n+d\leq 2d$$ avec $\rho= 2d$ si et seulement si $d$ est paire et $X=(S_3)^{d/2}$ o\`u $S_3$ est l'\'eclatement de $\Pbb^2$ en trois points non align\'es.
\end{teo}
Les preuves de ces deux r\'esultats sont inspir\'ees de celles des r\'esultats analogues dans le cas torique (\cite{De03}, \cite{Ca06}). Il faut prendre en compte le fait que les polytopes $G/H$-r\'eflexifs ne sont pas \`a sommets entiers comme dans le cas torique. Cependant, les sommets non entiers sont parmi un nombre fini de points rationnels qui d\'ependent seulement de $G/H$, ce qui permet de contr\^oler les changements qui apparaissent entre les cas torique et horosph\'erique. On utilisera alors des arguments de g\'eom\'etrie convexe ainsi que des \'el\'ements de combinatoire sur les groupes alg\'ebriques r\'eductifs.

On remarquera que les vari\'et\'es de nombre de Picard~$1$ sont souvent \'etudi\'ees \`a part: quelles sont ces vari\'et\'es? Les seules vari\'et\'es toriques lisses et de nombre de Picard $1$ sont les espaces projectifs. Par contre, parmi les vari\'et\'es horosph\'eriques lisses de nombre de Picard~$1$, on compte les vari\'et\'es de drapeaux $G/P$ avec $P$ maximal, mais aussi des vari\'et\'es non homog\`enes. L'\'etude de ces vari\'et\'es est en cours.\\

Dans la cinqui\`eme et derni\`ere partie, on d\'emontre le r\'esultat suivant.
\begin{teo}\label{tresample}
Soient $X$ une vari\'et\'e horosph\'erique projective de rang $n$, et $D$ un diviseur de Cartier et ample. Alors $(n-1)D$ est tr\`es ample.\\
Si de plus $X$ est localement factorielle, alors $D$ est tr\`es ample.
\end{teo}
Dans le cas torique, la premi\`ere partie de ce th\'eor\`eme est due \`a G.~Ewald et U.~Wessels \cite{EW91}, et la deuxi\`eme partie \`a M.~Demazure \cite{De70}. On utilisera un r\'esultat combinatoire de \cite{EW91} pour d\'emontrer ce r\'esultat.\\

Une question naturelle se pose: est-ce que ces trois th\'eor\`emes peuvent se g\'en\'eraliser aux vari\'et\'es sph\'eriques?

La premi\`ere assertion du th\'eor\`eme \ref{tresample} reste vraie pour les vari\'et\'es sph\'eri\-ques.
Quant au reste, on sait que toute vari\'et\'e sph\'erique de Fano $X$ d\'eg\'en\`ere en une vari\'et\'e horosph\'erique (ou m\^eme torique) $X_0$ qui est projective et $\Qbb$-Fano, c'est-\`a-dire, il existe un entier positif $k$ tel que $-kK_{X_0}$ est de Cartier et ample \cite{BA04}. Mais l'entier $k$ peut \^etre tr\`es grand, et la vari\'et\'e $X_0$ est en g\'en\'eral tr\`es singuli\`ere. Ainsi tous ces r\'esultats ne sont qu'une premi\`ere \'etape dans la classification des vari\'et\'es sph\'eriques de Fano.\\

Pour avoir des exemples de polytopes $G/H$-r\'eflexifs et de vari\'et\'es horosph\'eriques de Fano de rang~$2$, on peut regarder les chapitres 6 et 7 de \cite{Pa06}. On y trouve notamment une description des plongements lisses de Fano de $(SL_2\times\Cbb^*)/U$ et $SL_3/U$.

\section{Notations}\label{notations}
Toutes les vari\'et\'es consid\'er\'ees sont des vari\'et\'es alg\'ebriques sur $\Cbb$.

On se donne un groupe alg\'ebrique $G$ r\'eductif (c'est-\`a-dire qui ne contient aucun sous-groupe distingu\'e isomorphe \`a $\Cbb^n$) et connexe sur $\Cbb$, un sous-groupe de Borel $B$ de $G$, un tore maximal $T$ de $B$ et le radical unipotent $U$ de $B$. On note $R$ l'ensemble des racines de $(G,T)$, $R^+$ l'ensemble des racines positives (c'est-\`a-dire l'ensemble des racines de $(B,T)$), $S$ l'ensemble des racines simples, $\Lambda$ (respectivement $\Lambda^+$) le groupe des caract\`eres de $B$ ou de $T$ (respectivement l'ensemble des caract\`eres dominants) et  $W$ le groupe de Weyl de $(G,T)$. Pour toute racine simple $\alpha$, on note $\check\alpha$ sa coracine et $\omega_\alpha$ le poids fondamental associ\'e \`a $\alpha$.

Pour tout sous-groupe ferm\'e $H$ de $G$, $N_G(H)$ d\'esigne le normalisateur de $H$ dans $G$, et $R_u(H)$ est le radical unipotent de $H$.\\

Lorsque $I\subset S$, on note $W_I$ le sous-groupe de $W$ engendr\'e par les r\'eflexions simples $s_\alpha$ pour tout $\alpha\in I$, et de m\^eme $R_I$ (respectivement $R_I^+$) d\'esigne l'ensemble des racines (respectivement positives) qui sont combinaisons lin\'eaires des racines simples de $I$. On note $P_I$ le sous-groupe parabolique de $G$ engendr\'e par $B$ et $W_I$. Alors $I\longmapsto P_I$ est une bijection entre l'ensemble des parties de $S$ et l'ensemble des sous-groupes paraboliques contenant $B$ \cite[th. 8.4.3]{Sp98}.\\

Pour tout caract\`ere dominant $\lambda$, on note $V(\lambda)$ le $G$-module simple de plus grand poids $\lambda$ \cite[ch.XI]{Hu75}, et $V(\lambda)^*$ son dual. On  d\'esigne par $v_\lambda$ un vecteur propre de $V(\lambda)$ de poids $\lambda$, et le stabilisateur de la droite $\Cbb v_\lambda$ est un sous-groupe parabolique de $G$ contenant $B$ qu'on note $P(\lambda)$. En \'ecrivant $\lambda=\sum_{\alpha\in S}x_\alpha\omega_\alpha$, les $x_\alpha$ \'etant des entiers positifs ou nuls, on a $P(\lambda)=P_I$ o\`u $I$ est l'ensemble des racines simples $\alpha$ telles que $x_\alpha$ soit nul.

Les $G$-modules consid\'er\'es seront toujours rationnels et de dimension finie.
Si $V$ est un $G$-module, on note $V^G$ (respectivement $V^U$) l'ensemble des points fixes de $V$ sous l'action de $G$ (respectivement $U$), et $V^{(B)}$ l'ensemble des vecteurs propres de $V$ sous l'action de $B$. Comme $V$ est semi-simple, on a une d\'ecomposition en $G$-modules simples (o\`u les $m_\lambda$ sont des entiers positifs ou nuls): $$V=\bigoplus_{\lambda\in\Lambda^+}V(\lambda)^{\bigoplus m_\lambda}\quad\mbox{ et }\quad V^U=\bigoplus_{\lambda\in\Lambda^+}(\Cbb v_\lambda)^{\bigoplus m_\lambda}.$$
\begin{defi} Un sous-groupe ferm\'e $H$ de $G$ contenant $U$ est dit {\it horosph\'erique}. Dans ce cas, on dit aussi que l'espace homog\`ene $G/H$ est horosph\'erique.
\end{defi}
\begin{ex}
Soit $P$ un sous-groupe parabolique de $G$ contenant $B$ et soient $\chi_1,\dots,\chi_n$ des caract\`eres de $P$. Alors l'intersection des noyaux des $\chi_i$ dans $P$ est un sous-groupe horosph\'erique.
\end{ex}
En fait tout sous-groupe horosph\'erique $H$ est de cette forme.
\begin{prop}\label{horker}
Soit $H$ un sous-groupe horosph\'erique de $G$. Il existe un unique sous-groupe parabolique $P$ contenant $B$ tel que $H$ soit l'intersection de noyaux de caract\`eres de $P$. De plus, $P=N_G(H)$.
\end{prop}
\begin{proof}
D'apr\`es le th\'eor\`eme de Chevalley \cite[11.2]{Hu75}, il existe un $G$-module $V$ et une droite $L$ de $V$ telles que $H$ soit le stabilisateur de $L$, c'est-\`a-dire $H=\{x\in G\mid x.L=L\}$.

D\'ecomposons $V$ en somme directe de $G$-modules simples: $$V=\bigoplus_{\lambda\in\Lambda^+_0} V(\lambda)^{\oplus m_\lambda},$$ avec $m_\lambda$ non nul pour tout $\lambda$ dans un sous-ensemble fini $\Lambda^+_0$ de $\Lambda^+$.  Comme $U\subset H$, on a $L\subset V^U=\oplus_{\lambda\in\Lambda^+_0}(\Cbb v_\lambda)^{\oplus m_\lambda}$. Soit $\Lambda_1^+\subset\Lambda^+_0$ un sous-ensemble minimal tel que $L\subset \oplus_{\lambda\in\Lambda^+_1}(\Cbb v_\lambda)^{\oplus m_\lambda}$. Alors il existe $V'=\oplus_{\lambda\in\Lambda_1^+}V(\lambda)\subset V$ tel que la projection $L'$ de $L$ sur $V'$ v\'erifie $H=\{x\in G\mid x.L'=L'\}$. On peut donc supposer que $L$ est engendr\'ee par un vecteur de la forme $\sum_{\lambda\in\Lambda^+_1}a_\lambda v_\lambda$ avec tous les $a_\lambda\neq 0$. Notons $P=\bigcap_{\lambda\in\Lambda^+_1}P(\lambda)$.\\
Montrons que $P$ convient. Le th\'eor\`eme de Chevalley nous dit que
\begin{eqnarray*}
H & = & \{x\in G\mid\exists\lambda_0(x)\in \Cbb^*,\, x.v=\lambda_0(x)v\}\\
  & = & \{x\in\bigcap_{\lambda\in\Lambda^+_1}P(\lambda)\mid \forall\lambda,\mu\in\Lambda^+_1,\,\lambda(x)=\mu(x)\}.
\end{eqnarray*}
Soit $\mu\in \Lambda^+_1$, alors $H=\bigcap_{\lambda\in\Lambda^+_1}\ker(\lambda -\mu)\subset P$.

Montrons que $P=N_G(H)$. On a clairement $P \subset N_G(H)$, de plus $R_u(H)=R_u(P)$ et $P=N_G(R_u(P))$, donc $N_G(H)\subset N_G(R_u(H))=P$.
\end{proof}

\begin{defi} \label{INM} Soit $H$ un sous-groupe horosph\'erique.
On note $I$ le sous-ensemble de $S$ tel que $P=P_I$.
Puis on d\'efinit $M$  comme l'ensemble des caract\`eres de $P$ dont la restriction \`a $H$ est triviale; c'est un sous-r\'eseau de $\Lambda$. On note $N$ le r\'eseau dual de $M$.
Le rang de $M$ est appel\'e le rang de $G/H$; on le note $n$. On notera aussi $d$ la dimension de $G/H$, on a \'evidemment
\begin{equation}\label{dimension}
d=n+\dim(G/P)=n+\sharp(R^+\backslash R_I^+).
\end{equation}
On pose $M_\Rbb:=M\otimes_\Zbb\Rbb$ et $N_\Rbb:=N\otimes_\Zbb\Rbb$.
\end{defi}
\begin{rem}
L'espace homog\`ene $G/H$ est l'espace total d'une fibration sur la vari\'et\'e de drapeaux $G/P$ de fibre le tore $P/H$. Ce dernier est  isomorphe au tore dual de $M$ par l'application:
$$\begin{array}{ccc}
 P/H & \lra & \{\mbox{homomorphismes de groupes }M\rightarrow\Cbb^*\} \\
  pH & \lmt & [\chi \mapsto \chi(p)].
\end{array}$$
\end{rem}
On va classifier les espaces homog\`enes horosph\'eriques $G/H$ en termes de sous-r\'eseaux de $\Lambda$ et de sous-ensembles de $S$.
\begin{prop}
La constuction ci-dessus qui \`a un espace homog\`ene horosph\'erique $G/H$ associe le couple $(M,I)$ d\'efinit une bijection de l'ensemble des $G$-espaces homog\`enes horosph\'eriques sur l'ensemble des couples $(M,I)$, o\`u $I$ est un sous-ensemble de $S$ et $M$ un sous-r\'eseau de $\Lambda$ tel que pour tous $\alpha\in I$ et $\chi\in M$, $\langle\chi,\check\alpha\rangle=0$.
\end{prop}
\begin{proof}
Remarquons d'abord que pour tout $\chi\in\Lambda$, la condition $\langle\chi,\check\alpha\rangle=0$ pour $\alpha\in I$ est \'equivalente  au fait que $\chi$ s'\'etend en un caract\`ere de $P_I$.

Ensuite, \`a un couple $(M,I)$ v\'erifiant cette condition, on associe l'espace homog\`ene horosph\'eri\-que $G/H$, o\`u $H$ est l'intersection des noyaux des caract\`eres $\chi\in M$ dans $P_I$. On v\'erifie alors facilement que les deux applications sont inverses l'une de l'autre.
\end{proof}

\begin{defi}\label{plongiso}
Un {\it plongement} d'un espace homog\`ene $G/H$ est un couple $(X,x)$, o\`u $X$ est une $G$-vari\'et\'e alg\'ebrique normale et $x$ est un point de $X$, tels que  l'orbite de $x$ dans $X$ soit ouverte et isomorphe \`a $G/H$.\\
Deux plongements $(X,x)$ et $(X',x')$ sont isomorphes s'il existe un isomorphime $G$-\'equivariant de $X$ sur $X'$ qui envoie $x$ sur $x'$.

Une {\it vari\'et\'e horosph\'erique} est une $G$-vari\'et\'e alg\'ebrique normale qui con\-tient une orbite ouverte isomorphe \`a un espace homog\`ene horosph\'erique. Le rang d'une vari\'et\'e horosph\'erique est le rang de sa $G$-orbite ouverte.

Un espace homog\`ene $G/H$ est dit {\it sph\'erique} s'il contient une orbite ouverte sous l'action d'un sous-groupe de Borel $B$ de $G$. Une {\it vari\'et\'e sph\'erique} est une $G$-vari\'et\'e alg\'ebrique normale qui contient une orbite ouverte isomorphe \`a un espace homog\`ene sph\'erique. Toute vari\'et\'e horosph\'erique est sph\'erique; ceci r\'esulte en effet de la d\'ecomposition de Bruhat.
\end{defi}
\begin{rem}
Soit $(X,x)$ est un plongement d'un espace homog\`ene sph\'e\-rique (respectivement horosph\'erique) $G/H$. Alors $X$ est une vari\'et\'e sph\'erique (respectivement horosph\'erique).

Inversement, soit $X$ une vari\'et\'e sph\'erique (respectivement horosph\'erique).  Soit $x$ un point de l'orbite ouverte de $X$. Notons $H$ le stabilisateur de $x$ dans $G$. Alors, $G/H$ est un espace homog\`ene sph\'erique (respectivement horosph\'erique), et $(X,x)$ est un plongement de $G/H$. Il faut remarquer que la classe d'isomorphisme du plongement $(X,x)$ d\'epend du choix de $x$. En effet, soit $x'$ un autre point de l'orbite ouverte de $X$. Les plongements $(X,x)$ et $(X,x')$ sont  isomorphes si et seulement si le stabilisateur de $x'$ dans $G$ est aussi $H$; autrement dit si et seulement si on a $x'=p.x$, o\`u $p\in P=N_G(H)$.\\
Ainsi, il ne faut pas confondre \og classe d'isomorphisme de plongements de $G/H$ \fg~et \og classe d'isomorphisme de vari\'et\'es horosph\'eriques dont l'orbite ouverte est isomorphe \`a $G/H$\fg.

Dans la suite de l'article, le point d'un plongement est sous-entendu: \og soit $X$ un plongement de $G/H$\fg~signifie rigoureusement \og soit $(X,x)$ un plongement de $G/H$\fg.\\
\end{rem}
La classification des plongements d'un espace homog\`ene sph\'erique fix\'e est obtenue par l'\'etude de leurs orbites sous les actions de $G$ et $B$, mais aussi de leurs diviseurs irr\'eductibles stables sous ces actions.
\begin{defi}\label{couleurs}
Soit $G/H$ un espace homog\`ene sph\'erique. On note $\mathcal{D}$ l'ensemble des diviseurs irr\'eductibles de $G/H$ qui sont stables par $B$ mais non par $G$. Les \'el\'ements de $\mathcal{D}$ sont appel\'es {\it couleurs}.

Soit $X$ un plongement de $G/H$.
On note $X_1,\dots,X_m$ les diviseurs irr\'educ\-tibles de $X$ stables par $G$.
On peut identifier $\mathcal{D}$ avec l'ensemble des diviseurs irr\'eductibles de $X$ qui sont stables par $B$ mais non par $G$. Ainsi,  $\mathcal{D}\cup\{X_1,\dots,X_m\}$ est l'ensemble des diviseurs irr\'eductibles $B$-stables de $X$.

Une {\it couleur de $X$} est une couleur qui contient une $G$-orbite ferm\'ee.
\end{defi}
\begin{rem}
Attention, l'ensemble des couleurs de $G/H$ est vide. En effet, toute couleur $D\in\mathcal{D}$ est de codimension~$1$ donc $D$ ne contient pas $G/H$ (unique orbite ferm\'ee du plongement $G/H$). Lorsqu'on aura besoin de pr\'eciser \`a partir de quel espace homog\`ene horosph\'erique sont d\'efinies les couleurs, on dira: \og les couleurs associ\'ees \`a l'espace homog\`ene\fg.
\end{rem}
Si $G/H$ est horosph\'erique, l'ensemble des $B$-orbites de codimension $1$ de $G/H$ est l'ensemble des $Bw_0s_\alpha P/H$ lorsque $\alpha$ d\'ecrit $S\backslash I$ et o\`u $w_0$ est l'\'el\'ement de longueur maximale dans $W$. Les couleurs sont alors les adh\'eren\-ces $D_\alpha$ des $B$-orbites $Bw_0s_\alpha P/H$ dans $G/H$ et $\mathcal{D}$ est en bijection avec $S\backslash I$.
\begin{defi}\label{simpletoroidal}
Une vari\'et\'e sph\'erique est {\it simple} si elle ne contient qu'une seule orbite ferm\'ee. Si $G/H$ est un espace homog\`ene sph\'erique, alors tout plongement de $G/H$ est recouvert par les plongements simples de $G/H$ qu'il contient.

Une vari\'et\'e sph\'erique est {\it toro\"idale} si elle n'a aucune couleur.
\end{defi}
Soit $X$ un plongement d'un espace homog\`ene horosph\'erique $G/H$. D\'efi\-nissons une application
\begin{equation}\label{sigma}
\sigma: \mathcal{D}\cup\{X_1,\dots,X_m\}\lra N
\end{equation}
de la fa\c{c}on suivante\footnote{Elle est d\'efinie de la m\^eme fa\c{c}on dans le cas sph\'erique, o\`u on pose $M=\Cbb(G/H)^{(B)}/\Cbb^*$ (pour plus de d\'etails, se r\'ef\'erer \`a \cite{Kn91}).}. Soit $D$ un diviseur $B$-stable de $X$. Il d\'efinit naturellement une valuation $v_D$, $B$-invari\-ante, du corps des fonctions rationnelles $\Cbb(G/H)=\Cbb(X)$. On en d\'eduit donc un homomorphisme de groupes $\Cbb(G/H)^{(B)}/\Cbb^*\longrightarrow \Zbb$. En remarquant  ensuite que $M$ est isomorphe \`a \\$\Cbb(G/H)^{(B)}/\Cbb^*$, la restriction de $v_D$ \`a $\Cbb(G/H)^{(B)}/\Cbb^*$ d\'efinit alors un \'el\'ement de $N$ qu'on note $\sigma(D)$.
Notons que la restriction de $\sigma$ \`a $\mathcal{D}$ ne d\'epend pas de $X$ mais que de $G/H$. En fait, si $\alpha\in S\backslash I$, l'image par $\sigma$ de la couleur $D_\alpha$ est simplement la restriction \`a $M$ de $\check{\alpha}$. Dans ce cas, on notera cette image $\check\alpha_M$ au lieu de $\sigma(D_\alpha)$.
\begin{rem}
Il se peut que l'application $\sigma$ ne soit pas injective. Ainsi, si $G/H$ est horosph\'erique, $\sigma$ n'est pas toujours une bijection entre $\mathcal{D}$ et l'ensemble $\{\check\alpha_M\mid\alpha\in S\backslash I\}$. Par exemple, lorsque $H=P$, l'application $\sigma$ est constante car $N=\{0\}$.
\end{rem}
\begin{exs}\label{SL2/U}
(1) L'espace homog\`ene horosph\'erique $SL_2/U$, de rang $1$, est isomorphe \`a $\Cbb^2\backslash\{0\}$. On peut choisir $B$ (respectivement $U$) \'egal \`a l'ensemble des matrices de $SL_2$ de la forme $\left(\begin{array}{rr}
* & *\\
0 & *
\end{array}\right)$
 (respectivement $\left(\begin{array}{rr}
1 & *\\
0 & 1
\end{array}\right)$). Ici, $P=B$, $S=\{\alpha\}$, $I=\varnothing$ et $U$ est le noyau de $\omega_\alpha$ dans $B$. On remarque que le morphisme $SL_2/U\lra SL_2/P$ est la projection de $\Cbb^2\backslash\{0\}$ sur $\Pbb^1$.

L'action naturelle de $SL_2$ sur $\Cbb^2$ induit une action de $SL_2$ sur $\Pbb^2\simeq \Pbb(\Cbb\oplus\Cbb^2)$.  En notant $x_0,x_1,x_2$ les coordonn\'ees homog\`enes sur $\Pbb^2$, on remarque que $\Pbb^2$ est un plongement de $SL_2/U$. En fait $SL_2/U$ correspond \`a l'ouvert $\{[1,x_1,x_2],(x_1,x_2)\in\Cbb^2\backslash\{0\}\}$ de $\Pbb^2$. Notons encore $0$ le point fixe  $[1,0,0]$ de $\Pbb^2$ sous l'action de $SL_2$, $D$ la droite $\{[x_0,x_1,x_2]\in\Pbb^2, x_0=0\}$ (de sorte que $\Pbb^2\backslash D=\Cbb^2$) , et $E$ le diviseur exceptionnel de l'\'eclatement de $\Pbb^2$ au point $0$. Alors les plongements non triviaux de $SL_2/U$ sont les 5 vari\'et\'es pr\'esent\'ees dans le tableau suivant.
\begin{center}
\begin{tabular}{|c|c|c|c|c|}
\hline
&  plongement $X$ &  diviseur(s)  & $SL_2$-orbite(s) & couleur \\
&  de $SL_2/U $ & $SL_2$-stable(s) & ferm\'ee(s) & de $X$\\
\hline
\hline
1/ & $\Cbb^2$ & aucun & $\{0\}$ & $D_\alpha$\\
\hline
2/ & $\Pbb^2\backslash\{0\}$ & $D$  & $D$ & aucune \\
\hline
3/ & $\Pbb^2$ & $D$ & $D$ et $\{0\}$ & $D_\alpha$\\
\hline
4/ & $\Cbb^2$ \'eclat\'e en $0$ & $E$ & $E$ & aucune\\
\hline
5/& $\Pbb^2$ \'eclat\'e en $0$ & $D$ et $E$ & $D$ et $E$ & aucune\\
\hline
\end{tabular}
\end{center}

Pour l'espace homog\`ene $SL_2/U$,  $\mathcal{D}$ est le singleton $\{D_\alpha\}$ o\`u $D_\alpha$ est l'ensemble $\{[1,x_1,0],x_1\in\Cbb^*\}$.
Les plongements 1/, 2/ et 4/ n'ont qu'une $SL_2$-orbite ferm\'ee; ce sont des vari\'et\'es horosph\'eriques simples. Les plongements 2/, 4/ et 5/ n'ont pas de couleur; ce sont des vari\'et\'es horosph\'eriques toro\"idales (voir la d\'efinition \ref{simpletoroidal}). On remarque aussi que $\Pbb^2$ \'eclat\'e en $0$ est recouvert par $\Cbb^2$ \'eclat\'e en $0$ et $\Pbb^2\backslash\{0\}$, puis que $\Pbb^2$ est recouvert par $\Pbb^2\backslash\{0\}$ et $\Cbb^2$.

Dans cet exemple, les r\'eseaux $M$ et $N$ sont isomorphes \`a $\Zbb$. La figure suivante repr\'esente la droite $N_\Rbb$ avec l'image par $\sigma$ des diviseurs irr\'eductibles stables par $B$ de $\Pbb^2$ \'eclat\'e en $0$:
\begin{center}
\includegraphics[width=13cm]{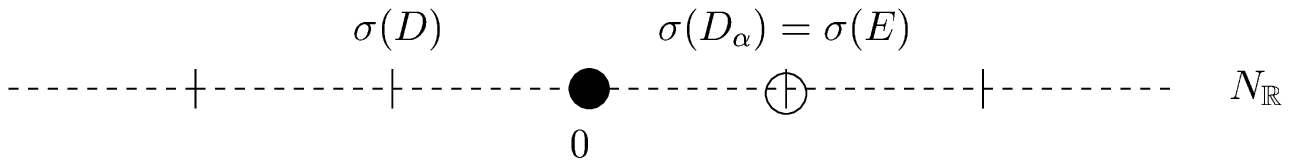}
\end{center}

(2) Quand $G/H$ est horosph\'erique, on rappelle que $P/H$ est isomorphe au tore dual de $M$ par l'application $p\in P\mapsto (\chi\in M\mapsto \chi(p))$. Soit $Y$ une vari\'et\'e torique sous l'action de ce tore; $P$ agit alors sur $Y$.

Soit $G\times^PY$ le quotient de $G\times Y$ par la relation d'\'equivalence  $(g,y)\sim(gp^{-1},p.y)$ pour tout $g\in G$, $p\in P$ et $y\in Y$. Alors $X=G\times^PY$ est une vari\'et\'e alg\'ebrique normale munie d'une fibration $G\times^P Y\longrightarrow G/P$. C'est aussi un plongement de $G/H$ et les diviseurs $X_i$ sont les $G\times^PY_i$ o\`u les $Y_i$ sont les diviseurs irr\'eductibles de $Y$ stables par le tore. Et pour tout $\alpha$ dans $S\backslash I$, $D_\alpha$ est $\overline{Bw_0s_\alpha P}\times^PY$. On remarque alors que chaque couleur ne contient aucune $G$-orbite de $X$: c'est une vari\'et\'e horosph\'erique toro\"idale. En fait, les vari\'et\'es horosph\'eriques toro\"idales sont toujours de la forme ci-dessus, c'est-\`a-dire des fibr\'es sur une vari\'et\'e de drapeaux, de fibre une vari\'et\'e torique (cela r\'esulte du th\'eor\`eme \ref{LunaVust} et de l'exemple \ref{SL2C}(3)).
\end{exs}

Lorsqu'on se donne un espace homog\`ene horosph\'erique $G/H$ de rang $n$, on lui associe un sous-groupe parabolique $P$, un ensemble $I\subset S$ et un r\'eseau $N\subset N_\Rbb$ de rang $n$ (voir la proposition \ref{horker} et la d\'efinition \ref{INM}). On a aussi l'ensemble des couleurs $\mathcal{D}$ (voir la d\'efinition \ref{couleurs}), et l'application $\sigma:\mathcal{D}\lra N$ (\ref{sigma}).

\begin{defi}
Soit $G/H$ un espace homog\`ene horosph\'erique\footnote{Si $G/H$ est sph\'erique, la d\'efinition des c\^ones et \'eventails colori\'es est quasiment identique (voir \cite[chap. 4]{Kn91}).} fix\'e (avec toutes les donn\'ees associ\'ees ci-dessus).

Un {\it c\^one colori\'e} de $N_\Rbb$ est un couple $(\mathcal{C}, \mathcal{F})$ o\`u $\mathcal{C}$ est un c\^one convexe de $N_\Rbb$ et $\mathcal{F}$ est un sous-ensemble de $\mathcal{D}$ appel\'e {\it l'ensemble des couleurs} du c\^one colori\'e, tel que\\
(i) $\mathcal{C}$ est engendr\'e par un nombre fini d'\'el\'ements du r\'eseau $N$ et contient $\sigma(\mathcal{F})$,\\
(ii) $\mathcal{C}$ est saillant (c'est-\`a-dire ne contient aucune droite) et
$\sigma(\mathcal{F})$ ne contient pas l'origine.

Une {\it face colori\'ee} d'un c\^one colori\'e $(\mathcal{C}, \mathcal{F})$ est un couple $(\mathcal{C'}, \mathcal{F'})$ o\`u $\mathcal{C'}$ est une face du c\^one $\mathcal{C}$ et $\mathcal{F'}$ est l'ensemble des \'el\'ements de $\mathcal{F}$ dont l'image par $\sigma$ est dans $\mathcal{C'}$.

Un {\it \'eventail colori\'e} de $N_\Rbb$ est un ensemble fini $\mathbb{F}$ de c\^ones colori\'es tel que\\
(i) toute face colori\'ee d'un c\^one colori\'e de $\mathbb{F}$ est dans $\mathbb{F}$,\\
(ii) pour tout \'el\'ement $u$ de $N_\Rbb$, il existe au plus un c\^one colori\'e $(\mathcal{C}, \mathcal{F})$ de $\mathbb{F}$ tel que $u$  soit dans l'int\'erieur relatif de $\mathcal{C}$.

Un \'eventail colori\'e $\mathbb{F}$ est dit {\it complet} si pour tout \'el\'ement $x$ de $N_\Rbb$, il existe un c\^one colori\'e $(\mathcal{C}, \mathcal{F})$ de $\mathbb{F}$ tel que $x$ soit dans $\mathcal{C}$.

On dira qu'un \'el\'ement $D$ de $\mathcal{D}$ est une {\it couleur} de $\mathbb{F}$ s'il existe un c\^one colori\'e de $\mathbb{F}$ dont $D$ est une couleur.
\end{defi}
\begin{rems}
Lorsque $G$ est un tore, l'ensemble des couleurs $\mathcal{D}$ est vide, et on retrouve la d\'efinition d'un \'eventail.

La derni\`ere condition dans la d\'efinition d'un \'eventail colori\'e implique que l'intersection de deux c\^ones colori\'es est une face colori\'ee commune.

Si deux couleurs ont la m\^eme image par $\sigma$, alors il se peut qu'un c\^one colori\'e ne poss\`ede qu'une des deux couleurs.
\end{rems}

\begin{ex}
Revenons \`a l'exemple \ref{SL2/U}(1). Les \'eventails colori\'es non triviaux de $N_\Rbb$ sont les suivants:
\begin{center}
\includegraphics[width=13cm]{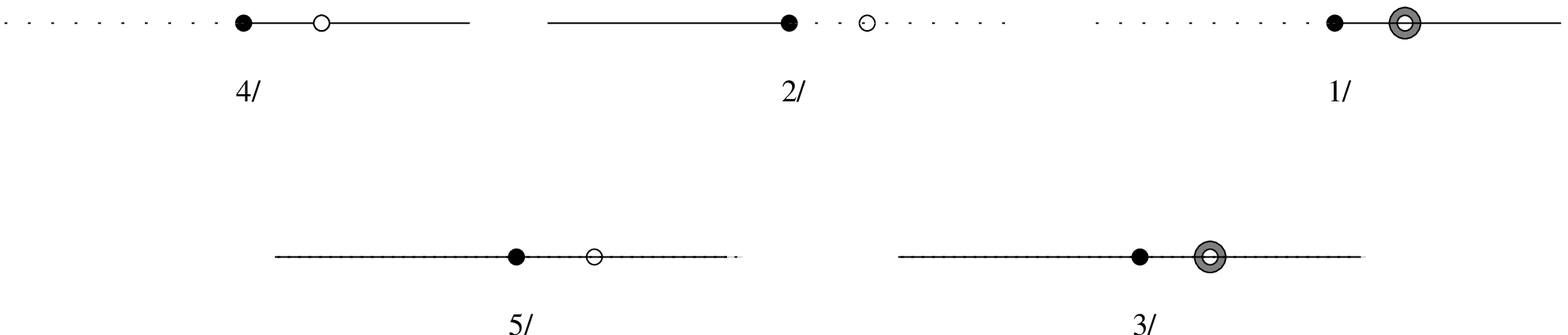}
\end{center}
On repr\'esente l'origine par un point noir, la couleur par un point blanc, et une couleur de l'\'eventail en ajoutant un anneau gris autour du point blanc. Les ar\^etes de l'\'eventail colori\'e sont les demi-droites noires issues de l'origine.

Comme $\sigma$ n'est pas injective, un point blanc pourrait \^etre l'image de deux couleurs. Dans toutes les figures de l'article, ce cas g\^enant n'aura pas lieu.\\
L'\'eventail i/ correspond au plongement i/ de $SL_2/U$ et l'\'eventail trivial $\{(\{0\},\varnothing)\}$ correspond au plongement trivial $SL_2/U$.
\end{ex}
On obtient un c\^one colori\'e $(\mathcal{C}, \mathcal{F})$ \`a partir d'un plongement simple $X$ de la fa\c{c}on suivante. Soit $Y$ l'unique $G$-orbite ferm\'ee de $X$; alors $\mathcal{F}$ est l'ensemble des couleurs qui contiennent $Y$. Ensuite $\mathcal{C}$ est le c\^one engendr\'e par $\sigma(\mathcal{F})$ et l'ensemble des $\sigma(D)$ lorsque $D$ parcourt l'ensemble des diviseurs irr\'eductibles de $X$ stables par $G$.

Pour un plongement quelconque $X$, l'\'eventail colori\'e associ\'e est l'ensemble des c\^ones colori\'es associ\'es aux plongements simples inclus dans $X$. On remarque que les couleurs de $X$ sont les m\^emes que les couleurs de l'\'eventail colori\'e associ\'e \`a $X$.
\begin{teo}[cas particulier\footnote{L'\'enonc\'e du th\'eor\`eme de D.~Luna et Th.~Vust \cite[th.4.3]{Kn91} est le m\^eme que celui du th\'eor\`eme ci-dessus avec $G/H$ sph\'erique.} du th\'eor\`eme 4.3 de \cite{Kn91}]
\label{LunaVust}
Soit $G/H$ un espace homog\`ene horosph\'erique.
La constuction ci-dessus d\'efinit une bijection entre l'ensemble des classes d'isomorphisme de plongements de $G/H$ (d\'efinition \ref{plongiso})  et l'ensemble des \'eventails colori\'es de $N_\Rbb$. De plus, les plongements complets correspondent aux \'eventails colori\'es complets.
\end{teo}
\begin{rems}
On retrouve la classification des vari\'et\'es toriques lorsque $G/H$ est un tore.
\end{rems}

\begin{exs}\label{SL2C}
(1) Voici un autre exemple de rang $1$. Notons $\alpha$ et $\beta$ les racines simples de $SL_3$ et posons $H=\ker(2\omega_\alpha-\omega_\beta)\subset B$. Les plongements complets de $SL_3/H$ sont en bijection  avec les \'eventails colori\'es suivants. \\
\begin{center}
\includegraphics[width=13cm]{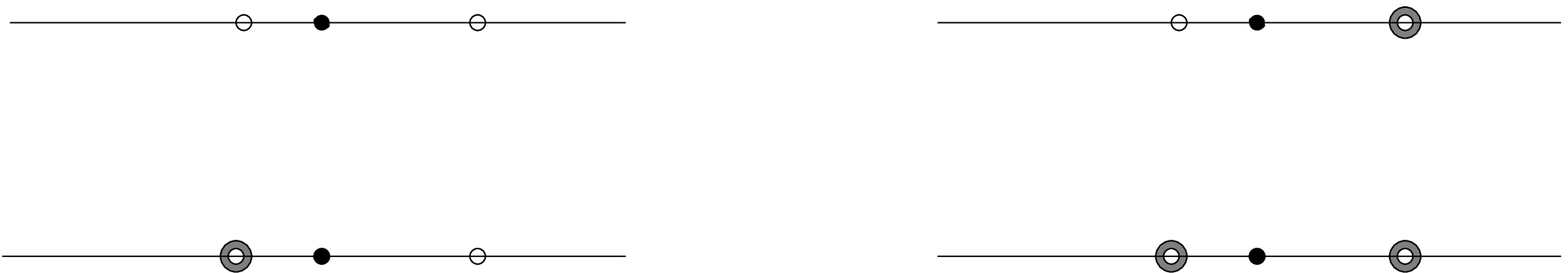}
\end{center}

(2) Donnons maintenant un exemple de rang $2$. Soient $G=SL_2\times\Cbb^*$ et $H=U$; alors $M$ a pour base $\omega_\alpha$ et le caract\`ere trivial de $\Cbb^*$, o\`u $\alpha$ est la racine simple de $SL_2$. Voici quelques exemples d'\'eventails colori\'es complets de $N_\Rbb$ qu'on peut obtenir dans ce cas.\\
\begin{center}
\includegraphics[width=13cm]{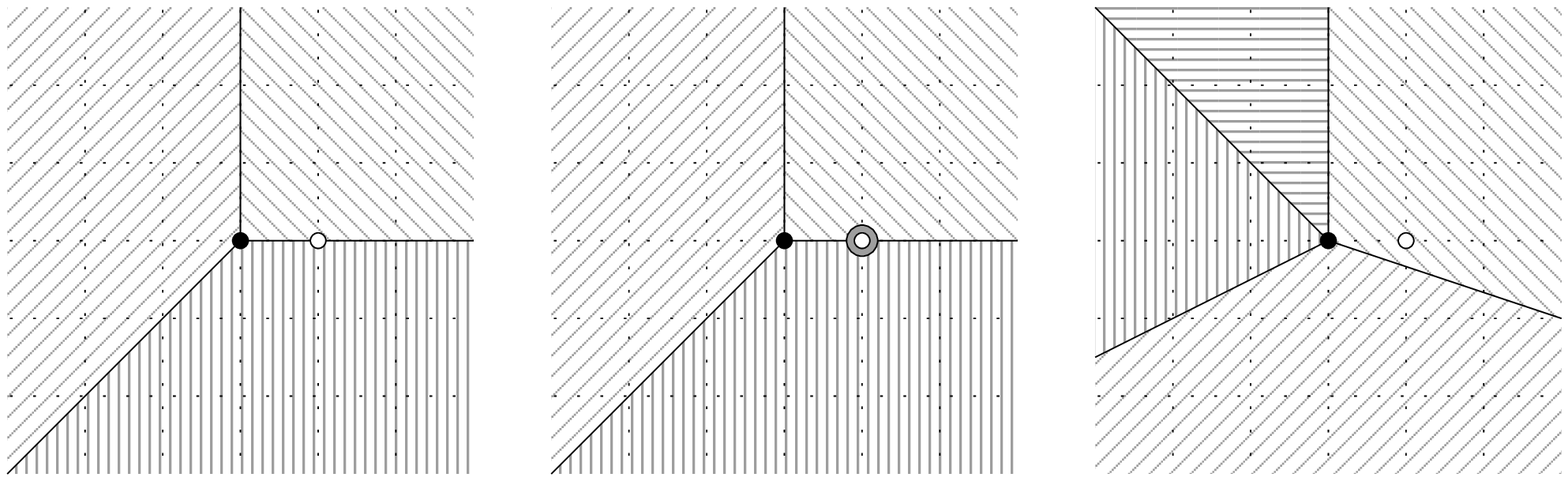}
\end{center}

On peut facilement remarquer que le nombre d'\'eventails colori\'es est fini lorsque le rang est $1$, alors qu'il est infini d\`es que le rang est au moins $2$.\\

(3) Soit $X=G\times^PY$ un plongement toro\"idal de $G/H$ comme dans l'exemple \ref{SL2/U} (2). Notons $\mathbb{E}$ l'\'eventail de $Y$ dans $N_\Rbb$; alors l'\'eventail colori\'e associ\'e \`a $X$ est l'ensemble des c\^ones colori\'es $(\mathcal{C},\varnothing)$ lorsque $\mathcal{C}$ parcourt $\mathbb{E}$. Dans l'exemple (1) ci-dessus, seul le premier \'eventail correspond \`a un plongement toro\"idal; ce dernier est un fibr\'e en $\Pbb^1$ sur $SL_3/B$. Et dans l'exemple (2) ci-dessus, le premier et le troisi\`eme \'eventails colori\'es correspondent \`a des plongements toro\"idaux.
\end{exs}

\section{Vari\'et\'es horosph\'eriques lisses}\label{sectionlisse}

F. Pauer a classifi\'e les plongements lisses de $G/H$ lorsque $H=U$ \cite{Pa83}. On va g\'en\'eraliser ce r\'esultat \`a tous les espaces homog\`enes horosph\'e\-riques $G/H$.
Mais avant de donner un crit\`ere de lissit\'e pour les vari\'et\'es horosph\'e\-riques, on va \'etudier une condition n\'ecessaire plus simple \`a caract\'eriser.
\begin{defi}
Une vari\'et\'e normale est dite {\it localement factorielle} si tout diviseur de Weil est de Cartier.
\end{defi}
Une vari\'et\'e lisse est toujours localement factorielle. La r\'eciproque est vraie dans le cas torique. Par contre elle ne l'est pas dans le cas sph\'erique, ni horosph\'erique. Par exemple, soit $\omega$ un poids fondamental, soit $X=\overline{G.v_\omega}$ le c\^one affine sur $G/P(\omega)$ dans $V(\omega)$; alors $X$ est un plongement de rang $1$ de $G/H$ pour $H=\ker(\omega)\subset P(\omega)$. On peut v\'erifier que $X$ est toujours localement factoriel; mais $X$ est lisse si et seulement si c'est $V(\omega)$ tout entier. On donnera d'autres exemples de vari\'et\'es horosph\'eriques localement factorielles et non lisses \`a la fin de cette partie.

A l'aide de la caract\'erisation des diviseurs de Cartier sur une vari\'et\'e sph\'erique \cite[prop.3.1]{Br89}, on va caract\'eriser les plongements d'un espace homog\`ene horosph\'erique $G/H$ qui sont localement factoriels. La preuve est laiss\'ee au lecteur.
\begin{prop}\label{locfac}
Soit $X$ un plongement d'un espace homog\`ene horosph\'erique $G/H$ d'\'eventail colori\'e $\mathbb{F}$. Alors $X$ est localement factoriel si et seulement pour tout c\^one colori\'e maximal $(\mathcal{C},\mathcal{F})$ de $\mathbb{F}$,\\
(i) les \'el\'ements de $\mathcal{F}$ ont des images deux \`a deux distinctes par $\sigma$,\\
(ii) $\mathcal{C}$ est engendr\'e par une partie d'une base de $N$ contenant $\sigma(\mathcal{F})$.
\end{prop}
\begin{rem}
Si $X=G\times^PY$ est un plongement toro\"idal (voir l'exemple \ref{SL2/U}(2)), l'\'eventail colori\'e $\mathbb{F}$ n'a pas de couleur. Donc $X$ est localement factoriel si et seulement si tout c\^one colori\'e maximal $(\mathcal{C},\varnothing)$ de $\mathbb{F}$ est engendr\'e par une base de $N$, c'est-\`a-dire si et seulement si $Y$ est localement factoriel (et m\^eme lisse). En fait, les plongements toro\"idaux localement factoriels d'un espace homog\`ene sph\'erique sont toujours lisses (ceci se d\'eduit de \cite[2.4 prop.1]{Br97b}).
\end{rem}
Avant d'\'enoncer la caract\'erisation des vari\'et\'es horosph\'eriques lisses, donnons la d\'efinition suivante.
\begin{defi}\label{lisse}
Soient $I$ et $J$ des sous-ensembles disjoints de $S$. Notons $\Gamma_S$ le diagramme de Dynkin de $G$, et $\Gamma_{I\cup J}$ le sous-graphe de $\Gamma_S$ dont les sommets sont les \'el\'ements de $I\cup J$ et les ar\^etes sont celles de $\Gamma_S$ qui relient deux \'el\'ements de $I\cup J$.

Alors on dira que $(I,J)$ est lisse si toute composante connexe $\Gamma$ de $\Gamma_{I\cup J}$ v\'erifie l'une des conditions suivantes:\\
1/ $\Gamma$ est un diagramme de Dynkin de type $A_n$ dont les sommets sont tous dans $I$ sauf une des deux extr\'emit\'es qui est dans $J$;
\begin{center}
\includegraphics{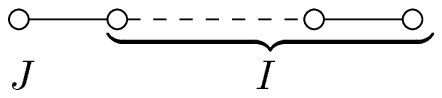}
\end{center}
2/ $\Gamma$ est un diagramme de Dynkin de type $C_n$ dont les sommets sont tous dans $I$ sauf l'extr\'emit\'e simple (c'est-\`a-dire non reli\'ee \`a l'ar\^ete double) qui est dans $J$;
\begin{center}
\includegraphics{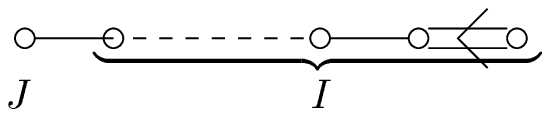}
\end{center}
3/ $\Gamma$ est un diagramme de Dynkin de type quelconque dont tous les sommets sont dans $I$.
\end{defi}
\begin{ex}
Consid\'erons le diagramme de Dynkin suivant:
\begin{center}
\includegraphics{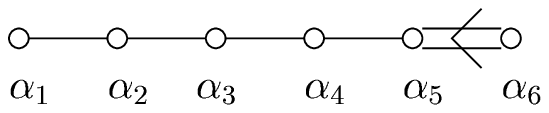}
\end{center}

Si $I=\{\alpha_2,\alpha_5,\alpha_6\}$ et $J=\{\alpha_1,\alpha_4\}$, alors $(I,J)$ est lisse.

Par contre, si $I=\{\alpha_1,\alpha_3\}$ et $J=\{\alpha_2\}$, alors $(I,J)$ n'est pas lisse.
\end{ex}

\begin{teo}\label{thlisse}
Soit $G/H$ un espace homog\`ene horosph\'erique. Rappelons que $S$ d\'esigne l'ensemble des racines simples et $I$ le sous-ensemble de $S$ associ\'e \`a $H$ (voir la d\'efinition \ref{INM}).

Un plongement $X$ de $G/H$, d'\'eventail colori\'e $\mathbb{F}$, est lisse si et seulement si les deux conditions suivantes sont v\'erifi\'ees.\\
(i) $X$ est localement factoriel (c'est-\`a-dire $\mathbb{F}$ v\'erifie les conditions de la proposition \ref{locfac}).\\
(ii) Pour tout c\^one colori\'e maximal $(\mathcal{C},\mathcal{F})$ de $\mathbb{F}$, on note $J_\mathcal{F}$ l'ensemble des $\alpha\in S\backslash I$ telles que $D_\alpha$ soit dans $\mathcal{F}$. Alors $(I,J_\mathcal{F})$ est lisse.
\end{teo}

\begin{rem}
Si $X$ est un plongement toro\"idal, il n'a pas de couleur. Or  $(I,\varnothing)$ est lisse pour tout $I$, donc la condition (ii) est toujours v\'erifi\'ee. On retrouve ainsi le fait que $X$ est lisse si et seulement si $X$ est localement factoriel.
\end{rem}

La preuve du th\'eor\`eme se fait en plusieurs \'etapes: apr\`es s'\^etre ramen\'e au cas o\`u $X$ est un plongement simple, on se ram\`enera au cas o\`u $X$ est affine, puis lorsque $X$ est lisse, au cas o\`u  c'est un $G$-module, et on fera ensuite une liste de tous les $G$-modules horosph\'eriques.

Avant de passer \`a la d\'emonstration du th\'eor\`eme, \'enon\c{c}ons quelques lemmes.

\begin{lem}\label{structlocale}
Soit $X$ un plongement simple de $G/H$. Notons $Y=G/H'$ son unique $G$-orbite ferm\'ee. On peut supposer que $H'\supset H$.
Alors il existe un unique morphisme $G$-\'equivariant $\pi:X\longrightarrow Y$.

Soient $P=N_G(H)$, $P'=N_G(H')$ et $L$ (respectivement $L'$) le sous-groupe de Levi de $P$ (respectivement $P'$) contenant $T$, si bien que $P\subset P'$ et $L\subset L'$.
Soit $Z=\pi^{-1}(P'/H')$; c'est une $L'$-vari\'et\'e. Alors $Z$ est une vari\'et\'e horosph\'erique sous l'action de $L'$, c'est un plongement affine (et donc simple) de $L'/(L'\cap H)$, dont l'unique $L'$-orbite ferm\'ee est le tore $P'/H'$. De plus le c\^one colori\'e de $Z$ s'identifie \`a celui de $X$.
\end{lem}
\begin{rem}
Si $X$ est une $G$-vari\'et\'e affine alors les $G$-orbites ferm\'ees de $X$ sont s\'epar\'ees par $\Cbb[X]^G$ (voir \cite[II.3.3]{Kr85}). Si de plus, $X$ est une vari\'et\'e sph\'erique, alors toute fonction de $\Cbb[X]^G$ est constante. Par cons\'equent, toute vari\'et\'e sph\'erique affine est simple.
\end{rem}
\begin{proof}
Construisons un morphisme $G$-\'equivariant $\pi:X\longrightarrow Y$.

On sait que $X$ est quasi-projectif; autrement dit, il existe un $G$-module $V$ tel que $X$ soit une sous-vari\'et\'e localement ferm\'ee de $\Pbb(V)$.
$$\xymatrix{
    G/H\, \ar@{^{(}->}[r]  & \,X\,   \ar@{^{(}->}[r]^\iota & \,\Pbb(V)\\
     & Y=G/H'  \ar@{^{(}->}[u] &
  }$$

Soient $v_0\in V$ tel que $[v_0]=\iota(y_0)$, $V'$ le sous-$G$-module de $V$ engendr\'e par $v_0$, et $\phi:\Pbb(V)\dashrightarrow\Pbb(V')$ l'application rationnelle $G$-\'equivariante d\'efinie par une projection $pr:V\longrightarrow V'$.

Soit $x\in X$; comme $X$ est un plongement simple, l'adh\'erence $\overline{G.x}$ dans $X$ de l'orbite de $x$ contient $Y$. Soit $v\in V$ tel que  $[v]=\iota(x)$, alors $v$ ne peut pas \^etre dans le noyau de $pr$. Par cons\'equent, l'application rationnelle $\phi\circ\iota$ est d\'efinie en $x$.

De plus, la droite $\Cbb v_0$ est stabilis\'ee par le sous-groupe $H'$ contenant $U$, donc $v_0$ est une somme de vecteurs propres de plus haut poids. Du fait que $ v_0\in\overline{G.v}$, on en d\'eduit que $\phi\circ\iota(x)\in G.[v_0]\subset\Pbb(V')$. Via l'isomorphisme $G.[v_0]\simeq G/H'$, on peut alors d\'efinir $\pi:X\longrightarrow G/H'$ par $\phi\circ\iota$.
$$\xymatrix{
X\, \ar@{^{(}->}[r]^\iota \ar@{->}[rd]^{\phi\circ\iota} \ar@/^1pc/[d]^\pi & \,\Pbb(V) \ar@{-->}[rd]^\phi & \\
G/H'\, \ar[r]^{\sim} \ar@{^{(}->}[u] & G.[v_0] \ar@{^{(}->}[r] & \,\Pbb(V')
}$$
La premi\`ere partie du lemme est alors montr\'ee.

On d\'efinit alors $Z$ comme dans l'\'enonc\'e. C'est naturellement une $L'$-vari\'et\'e normale car $\pi: X\longrightarrow G/P'$ est une fibration de fibre $Z$. L'orbite ouverte de $Z$ est $Z\cap (G/H)=P'/H=L'/(L'\cap H)$. La derni\`ere \'egalit\'e est donn\'ee par le fait que $R_u(P')\subset H\subset P'$  implique $H=(L'\cap H)R_u(P')$. De plus l'unique $L'$-orbite ferm\'ee de $Z$ est $P'/H'$, donc $Z$ est un plongement simple de $L'/(L'\cap H)$.

Pour montrer que $Z$ est affine, \'etudions les couleurs associ\'ees \`a $L'/(L'\cap H)$ et celles de $Z$.\\
Les couleurs associ\'ees \`a $L'/(L'\cap H)$ sont  les $D_\alpha'=D_\alpha\cap (L'/(L'\cap H))$ lorsque $\alpha$ d\'ecrit le sous-ensemble $J$ de $S\backslash I$ tel que $P'=P_{I\cup J}$.

Soit $\alpha\in J$; alors $\pi(D_\alpha)$ est l'adh\'erence dans $G/H'$ de $Bw_0s_\alpha PH'/H'=Bw_0s_\alpha P'/H'=Bw_0P'/H'$, c'est-\`a-dire $\pi(D_\alpha)=G/H'$. Donc l'adh\'erence de $D_\alpha$ dans $X$ contient $Y$, c'est-\`a-dire $D_\alpha$ est une couleur de $X$.

Au contraire, si $\alpha\in S\backslash (I\cup J)$ alors $\pi(D_\alpha)$ est de codimension $1$ dans $G/H'$, donc $D_\alpha$ n'est pas une couleur de $X$.

En r\'esum\'e, $J$ est l'ensemble des racines simples associ\'ees \`a une couleur de $X$. On a aussi montr\'e que toutes les couleurs associ\'ees \`a $L'/(L'\cap H)$ sont aussi des couleurs de $Z$, car pour tout $\alpha\in J$, l'adh\'erence de $D_\alpha'$ dans $Z$ contient $P'/H'$. Or un plongement simple dont les couleurs sont exactement celles associ\'ees \`a l'espace homog\`ene est affine \cite[th.3.1]{Kn91}, donc $Z$ est affine.

Il reste \`a montrer que  le c\^one colori\'e de $Z$ s'identifie \`a celui de $X$.

Rappelons que le r\'eseau $M$ associ\'e \`a $G/H$ est l'ensemble des caract\`eres $\chi$ de $P$ qui s'annulent sur $H$. Notons $M'$ le r\'eseau associ\'e \`a $L'/(L'\cap H)$, c'est-\`a-dire l'ensemble des caract\`eres de $L'\cap P$ qui s'annulent sur $L'\cap H$. On a imm\'ediatement $M\subset M'$. Montrons l'inclusion oppos\'ee: soit $\chi\in M'$, alors $\chi$ est un caract\`ere de $L'\cap P$ donc aussi de $(L'\cap P)R_u(P')=P$, de plus $\chi$ s'annule sur $L'\cap H$ donc aussi sur $(L'\cap H)R_u(P')=H$. On a donc $M=M'$.

Les c\^ones colori\'es de $X$ et $Z$ sont donc tous les deux dans le m\^eme espace $N_\Rbb$ et les images des couleurs $D_\alpha$ pour tout $\alpha\in J$, sont les m\^emes dans les deux cas.  On a d\'ej\`a vu que les couleurs de $X$ et $Z$ sont les m\^emes, il suffit donc de montrer qu'ils ont le m\^eme c\^one, ou encore que les diviseurs irr\'eductibles $G$-stables de $X$ sont les m\^emes que les diviseurs irr\'eductibles $L'$-stables de~$Z$. C'est bien le cas, car $\pi$ d\'efinit une fibration $G$-\'equivariante $X\longrightarrow G/P'$ de fibre $Z$.
\end{proof}
Comme le morphisme $\pi$ est une fibration de fibre $Z$ sur la vari\'et\'e de drapeaux $G/P'$, $X$ est lisse si et seulement si $Z$ est lisse.
Le lemme pr\'ec\'edent nous ram\`ene donc \`a une vari\'et\'e horosph\'erique affine dont la $G$-orbite ferm\'ee est un tore. On va maintenant se ramener \`a l'\'etude des $G$-modules horosph\'eriques.
\begin{lem}\label{pointfixe}
Soit $X$ un plongement affine (et donc simple) de $G/H$ dont l'unique $G$-orbite ferm\'ee $G/H'$ est un tore. Notons $H'^0$ la composante neutre de $H'$; c'est un sous-groupe r\'eductif connexe de $G$ contenant la partie semi-simple de $G$.

Alors $X$ est une fibration homog\`ene sur $G/H'$ de fibre une $H'^0$-vari\'et\'e horosph\'erique $X_0$ avec un point fixe.

De plus $X$ est lisse si et seulement si $X_0$ est isomorphe \`a un $H'^0$-module.
\end{lem}
\begin{proof}
Soit $\pi: X\lra G/H'$ l'application d\'efinie au lemme pr\'ec\'e\-dent. Notons $X_0$ la fibre de $\pi$ au dessus du point $H'/H'$. C'est une $H'^0$-vari\'et\'e horosph\'erique et on a $X\simeq G\times^{H'}X_0$.

De plus $X_0$ est affine car il est ferm\'e dans $X$. Le point $H'/H'$ est un point fixe de $X_0$ sous l'action de $H'^0$.

Il existe un $H'^0$-module $V$ tel que $X_0$ soit isomorphe (comme $H'^0$-vari\'et\'e) \`a un ferm\'e de $V$. Quitte \`a faire une translation, on peut supposer que le point fixe de $X_0$ est $0$. De plus $X_0$ est horosph\'erique, en particulier $X_0=H'^0.X_0^U$, donc $X_0$ est stable par action lin\'eaire de $\Cbb^*$ dans $V$ \`a poids positif. On peut supposer que $X_0$ engendre le $H'^0$-module $V$, l'espace tangent \`a $X_0$ en $0$ est alors \'egal \`a $V$.

Si $X$ est lisse, alors $X_0$ l'est aussi et la dimension de $V$ est \'egale \`a celle de $X_0$. On en d\'eduit donc que $X_0$ est le $H'^0$-module $V$.
\end{proof}
Le lemme suivant donne une condition n\'ecessaire pour qu'un $G$-module soit horosph\'erique, lorsque $G$ est semi-simple.
\begin{lem} \label{exobourbaki}
Supposons $G$ semi-simple. On note $G_1,\dots,G_k$ les sous-groupes simples distingu\'es de $G$. On sait alors que $G$ est le quotient du produit des $G_i$ par un sous-groupe fini central.

Soit $V$ un $G$-module horosph\'erique. Alors, quitte \`a r\'earranger les indices, le $G$-module $V$ est isomorphe \`a $\oplus_{i=1}^{k'} V(\lambda_i)$ o\`u $V(\lambda_i)$ est un $G_i$-module simple pour tout $i\in \{1,\dots,k'\}$, et $k'\leq k$.
\end{lem}
Ce lemme d\'ecoule directement de \cite[chap.7 ex. 18 ]{Bo75} et de la remarque suivante.
\begin{rem}
Soit $V=\oplus_{i=1}^{k'} V(\lambda_i)$ un $G$-module horosph\'erique. Alors pour tout $i\neq j$, on a $V(\lambda_i)^*.V(\lambda_j)^*=V(\lambda_i+\lambda_j)^*$ dans $\Cbb[V]$ vue comme l'alg\`ebre sym\'etrique de $V^*$.

En effet, comme $V$ est horosph\'erique, $\Cbb[V]$ est une sous-alg\`ebre de $\Cbb[G/U]$. De plus, comme $T$ normalise $U$, il agit par multiplication \`a droite sur $\Cbb[G/U]$. De plus on sait que pour tout $\lambda\in\Lambda^+$, $$V(\lambda)^*\simeq\{f\in\Cbb[G]\mid\forall g\in G,\,\forall b\in B,\, f(gb)=\lambda(b)f(g)\}\subset\Cbb[G/U].$$  On a ainsi un isomorphisme de $G\times T$-modules $\Cbb[G/U]\simeq\oplus_{\lambda\in\Lambda^+}V(\lambda)^*$ et $V(\lambda)^*.V(\mu)^*=V(\lambda+\mu)^*$ pour tous $\lambda$, $\mu\in\Lambda^+$.
\end{rem}
Le lemme pr\'ec\'edent nous ram\`ene au cas o\`u $G$ et $V$ sont simples, cas d\'ecrit dans le lemme suivant.
\begin{lem}\label{gmodhorosph}
Supposons que $G$ est simple. Un $G$-module simple $V(\lambda)$ est horosph\'erique si et seulement si on est dans l'un des deux cas suivants:
\begin{enumerate}
\item $G=SL_n$ c'est-\`a-dire le syst\`eme de racines est de type $A_{n-1}$, et $\lambda$ est le poids fondamental $\omega_\alpha$ o\`u $\alpha$ est l'une des deux racines simples situ\'ees aux extr\'emit\'es dans le diagramme de Dynkin;
\begin{center}\includegraphics{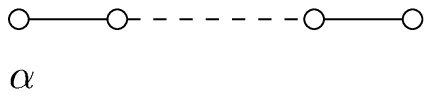}
\end{center}
\item $G=Sp_{2n}$ c'est-\`a-dire le syst\`eme de racines est de type $C_{n}$, et $\lambda$ est le poids fondamental $\omega_\alpha$ o\`u $\alpha$ est la racine simple situ\'ees \`a l'extr\'emit\'e simple dans le diagramme de Dynkin.
\begin{center}\includegraphics{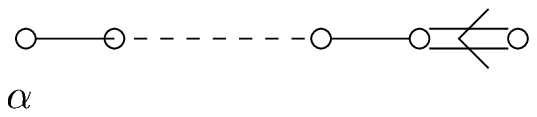}
\end{center}
\end{enumerate}
\end{lem}
\begin{proof}
Si $V(\lambda)$ est horosph\'erique alors  $G/P(\lambda)$ est le projectivis\'e de $V(\lambda)$, donc le groupe de Picard de $G/P(\lambda)$ est isomorphe \`a $\Zbb$ engendr\'e par $\mathcal{O}(1)$. Donc $P(\lambda)$ est un sous-groupe parabolique maximal de $G$, et $\lambda$ est un poids fondamental.

Les points fixes dans $G/P(\lambda)$ d'un tore maximal $T$ de $G$ sont tous conju\-gu\'es
par le groupe de Weyl $W$, donc les poids de $V(\lambda)$ le sont aussi :
ce sont les $w(\lambda)$, $w$ dans $W$. Cela signifie que $\lambda$ est un
poids minuscule au sens de \cite[chap 8, 7.3.]{Bo75}.

Il suffit alors de chercher, parmi la liste de ces cas \cite[chap 8, 7.3.]{Bo75}, ceux qui correspondent effectivement \`a des $G$-modules horosph\'eriques. Pour cela, on calcule les dimensions de $G/P(\lambda)$ et de $V(\lambda)$. Si celle de $G/P(\lambda)$ vaut celle de $V(\lambda)$ moins un, alors $V(\lambda)$ est horosph\'erique, et r\'eciproquement. On obtient, apr\`es calcul, le r\'esultat voulu.
\end{proof}
Il est maintenant possible de classifier les $G$-modules horosph\'eriques, lorsque $G$ est un groupe r\'eductif quelconque.
\begin{cor}[des deux lemmes pr\'ec\'edents]\label{cor2lem}
Soient $C$ la composante neutre du centre de $G$ et $G'$ la partie semi-simple de $G$. On rappelle que $G=C.G'$. On note aussi $G_1,\dots,G_k$ les sous-groupes simples distingu\'es de $G'$.

Soit $V$ un $G$-module; alors $V$ est horosph\'erique si et seulement si $V=\oplus_{i=1}^nV(\lambda_i)$ de sorte que\\
(1) $(\lambda_1,\dots,\lambda_n)$ est une famille libre de caract\`eres dominants,\\
(2) il existe $n'\leq n$ tel que, pour tout $i\in\{1,\dots,n'\}$, $V(\lambda_i)$ est un $G_i$-module simple comme dans le lemme \ref{gmodhorosph}, et pour tout $i\in\{n'+1,\dots,n\}$, $\lambda_i$ est un caract\`ere de $C$ (c'est-\`a-dire $V(\lambda_i)=\Cbb$).
\end{cor}
\begin{proof}
D\'ecomposons $V$ de la fa\c{c}on suivante: $V=V'\oplus V^{G'}$. On peut alors \'ecrire $V$ sous la forme $\oplus_{i=1}^nV(\lambda_i)$ de sorte que $V'=\oplus_{i=1}^{n'}V(\lambda_i)$ et $V^{G'}=\oplus_{i=n'+1}^nV(\lambda_i)$.

Ainsi, pour $i\in\{1,\dots,n'\}$, $V(\lambda_i)$ est un $G'$-module simple, et pour $i\in\{n'+1,\dots,n\}$, $V(\lambda_i)$ est un $C$-module simple ($V(\lambda_i)=\Cbb$).

Si $V$ est horosph\'erique alors $V'$ l'est sous l'action de $G'$ et $V^{G'}$ l'est sous l'action de $C$. Donc, par le lemme \ref{exobourbaki}, pour tout $i\in\{1,\dots,n'\}$, $V(\lambda_i)$ est un $G_i$-module simple (quitte \`a r\'earranger les indices) qui doit bien s\^ur v\'erifier le lemme \ref{gmodhorosph}. De plus, $(\lambda_{n'+1},\dots,\lambda_n)$ est une famille libre donc $(\lambda_1,\dots,\lambda_n)$ l'est aussi.

La r\'eciproque se d\'eduit du lemme \ref{gmodhorosph}.
\end{proof}

\begin{rem}\label{remgmodhorosph}
Si $V$ est comme dans le lemme ci-dessus, on note $J$ l'ensemble des racines simples $\alpha$ de $G$ telles que $\omega_\alpha$ soit l'un des $\lambda_i$.

Alors $V$ est un plongement simple de $G/H$ o\`u $H=\cap_{i=1}^n\ker\lambda_i\subset P=\cap_{i=1}^nP(\lambda_i)$. Le r\'eseau $M$ associ\'e \`a $V$ est le sous-r\'eseau de $\Lambda$ engendr\'e par la famille $(\lambda_1,\dots,\lambda_n)$, et le c\^one colori\'e $(\mathcal{C},\mathcal{F})$ associ\'e \`a $V$ est tel que $\mathcal{F}=\{D_\alpha,\, \alpha\in J\}$ et $\mathcal{C}$ est le c\^one engendr\'e par la base duale de $(\lambda_1,\dots,\lambda_n)$ dans $N$.

De plus, $(I,J)$ est lisse d'apr\`es le lemme \ref{gmodhorosph}.
\end{rem}

On a maintenant tout les outils pour achever la preuve du th\'eor\`eme \ref{thlisse}.\\
Comme tout plongement de $G/H$ est recouvert par des plongements simples et que son \'eventail colori\'e est donn\'e par l'ensemble des faces des c\^ones colori\'es de ses plongements simples, on en d\'eduit qu'il suffit de montrer le th\'eor\`eme dans le cas o\`u $X$ est un plongement simple.

On cherche alors les plongements simples de $G/H$ qui sont lisses. D'apr\`es les lemmes \ref{structlocale} et \ref{pointfixe}, cela revient \`a chercher les $H'^0$-modules horosph\'eriques o\`u $H'^0$ est un groupe r\'eductif dont le diagramme de Dynkin est $\Gamma_{I\cup J_\mathcal{F}}$.

Le corollaire \ref{cor2lem} montre l'implication ($X$ lisse $\Longrightarrow$ (ii)) du th\'eor\`eme \ref{thlisse}. L'implication ($X$ lisse$\Longrightarrow$ (i)) \'etant triviale, il suffit de montrer que (i) et (ii) impliquent que $X_0$ est un $H'^0$-module. Le c\^one colori\'e de $X_0$ est le m\^eme que celui de $X$. Le corollaire \ref{cor2lem} et la remarque \ref{remgmodhorosph} nous disent alors que $X_0$ a le m\^eme c\^one colori\'e qu'un $H'^0$-module; $X_0$ est donc isomorphe \`a ce $H'^0$-module.

\begin{exs}\label{exSL3U}
(1) Les deux plongements projectifs de $SL_2/U$ sont lisses. En revanche, le seul plongement projectif lisse de $SL_2/\ker(2\omega_\alpha)$ est le plongement toro\"idal; c'est m\^eme le seul plongement localement factoriel.

(2) Les plongements lisses de $(SL_2\times\Cbb)/U$ sont exactement les plongements localement factoriels.

(3) Regardons les plongements de $(SL_2\times SL_2)/U$. Le r\'eseau $N$ est de rang $2$ engendr\'e par les couleurs $\check\alpha_M=\check\alpha$ et  $\check\beta_M=\check\beta$ (o\`u $\alpha$ et $\beta$ sont les racines simples de $SL_2\times SL_2$). Ainsi les plongements lisses de $SL_2\times SL_2$ sont encore les plongements localement factoriels.

(4) Int\'eressons-nous maintenant aux plongements de $SL_3/U$. On obtient exactement les m\^emes \'eventails colori\'es que pour $(SL_2\times SL_2)/U$ car si $\alpha$ et $\beta$ sont les racines simples de $SL_3$, le r\'eseau $N$ est toujours engendr\'e par les couleurs $\check\alpha_M=\check\alpha$ et  $\check\beta_M=\check\beta$. Par contre les plongements localement factoriels ne sont pas toujours lisses: la condition d'\^etre localement factoriel ne d\'epend que de $N$ et de l'emplacement des couleurs dans $N$, par contre la lissit\'e d\'epend aussi des racines de $G$. L'\'eventail colori\'e suivant donne un plongement lisse de $(SL_2\times SL_2)/U$, mais aussi un plongement localement factoriel et non lisse de $SL_3/U$, car $(\varnothing, \{\alpha,\beta\})$ n'est pas lisse lorsque $\alpha$ et $\beta$ sont les racines simples de $SL_3$.
\begin{center}
\includegraphics[width=13cm]{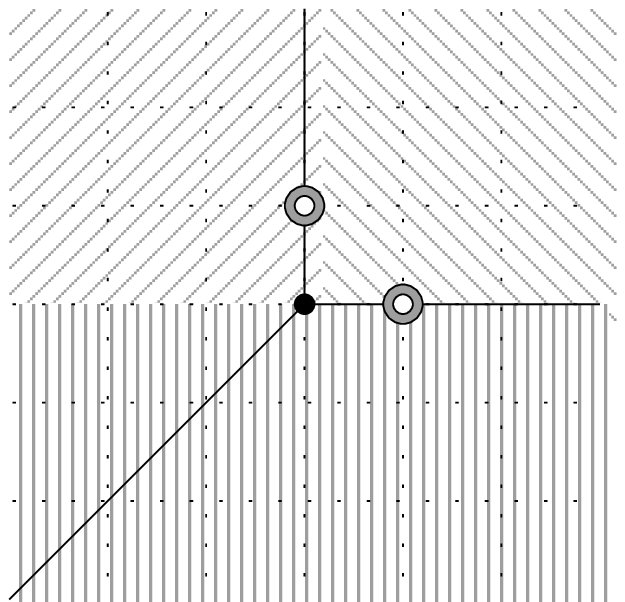}
\end{center}
Il faut aussi remarquer que, bien que les \'eventails colori\'es soient les m\^emes pour $(SL_2\times SL_2)/U$ et $SL_3/U$, ils repr\'esentent dans le premier cas des vari\'et\'es de dimension $4$, et dans le second, des vari\'et\'es de dimension $5$.
\end{exs}

A l'aide de la d\'emonstration du th\'eor\`eme \ref{thlisse}, on va montrer le r\'esultat suivant, d\'ej\`a connu dans le cas torique.
\begin{prop}\label{sousvarlisse}
Toute sous-vari\'et\'e irr\'eductible et stable par $G$ d'une vari\'et\'e horosph\'erique lisse est aussi lisse.
\end{prop}
\begin{proof}
Soient $X$ un plongement lisse d'un espace homog\`ene horosph\'erique $G/H$, et $Y$ une sous-vari\'et\'e irr\'eductible et stable par $G$ de $X$.

Si $Y$ est singuli\`ere, le lieu singulier de $Y$ est ferm\'e et stable par $G$, donc il contient une $G$-orbite ferm\'ee de $X$. Par cons\'equent il suffit de montrer que $Y$ est lisse le long de toute $G$-orbite ferm\'ee. On peut alors supposer que $X$ est un plongement simple de $G/H$.

Gr\^ace aux lemmes \ref{structlocale} et \ref{pointfixe}, on se ram\`ene au cas o\`u $X$ est un $H'^0$-module horosph\'erique et $Y$ une sous-vari\'et\'e irr\'eductible et stable par $H'^0$ (o\`u $H'^0$ est le groupe r\'eductif connexe d\'efini au lemme \ref{pointfixe}). Il suffit donc de montrer que les sous-vari\'et\'es irr\'eductibles et stables par $G$ d'un $G$-module horosph\'erique sont lisses (pour tout groupe r\'eductif connexe $G$).

Pour conclure, on va d\'ecrire les $G$-orbites des $G$-modules horosph\'eriques, en utilisant la description des $G$-modules horosph\'eriques donn\'ee dans le corollaire \ref{cor2lem}.

Avec les notations du corollaire \ref{cor2lem}, on sait que pour $i\in\{1,\dots,n'\}$, les $G_i$-orbites (ou $G$-orbites) de $V(\lambda_i)$ sont le point fixe $0$ et $V(\lambda_i)\backslash\{0\}$. Et pour $i\in\{n'+1,\dots,n\}$, les $G$-orbites de $V(\lambda_i)$ sont $0$ et $\Cbb\backslash\{0\}$. De plus, les caract\`eres $\lambda_i$ sont ind\'ependants, donc les $G$-orbites de $V$ sont les sommes partielles des $V(\lambda_i)\backslash\{0\}$.

On en d\'eduit que les sous-vari\'et\'es irr\'eductibles et stables par $G$ d'un $G$-module horosph\'erique $V$ sont les sous-$G$-modules de $V$. Elles sont donc \'evidemment lisses, ce qui montre la proposition.
\end{proof}
\begin{rem}
Ce r\'esultat n'est pas vrai pour les vari\'et\'es sph\'eriques (cf \cite{Br94}).
\end{rem}

\section{Classification des plongements de Fano}\label{classificationFano}
Le but de cette partie est de classifier, en termes de polytopes, les plongements projectifs d'un espace homog\`ene horosph\'erique fix\'e qui sont de Fano, de fa\c{c}on \`a g\'en\'eraliser la classification des vari\'et\'es toriques de Fano.

Rappelons qu'une vari\'et\'e projective est de Fano si elle est normale et si son diviseur anticanonique est de Cartier et ample. Dans cette partie, on fixe un espace homog\`ene horosph\'erique $G/H$. Les notations sont celles de la partie \ref{notations}.

Lorsque $G/H$ est de rang~$0$, c'est-\`a-dire lorsque $H$ est un sous-groupe parabolique $P$, le seul plongement est la vari\'et\'e de drapeaux $G/P$, qui est lisse et de Fano. Je me place maintenant dans le cas o\`u $G/H$ est de rang au moins $1$. Les \'eventails colori\'es complets consid\'er\'es ne seront donc pas r\'eduits au point~$0$.

La premi\`ere \'etape consiste \`a d\'eterminer un diviseur anticanonique des plongements de $G/H$.

\begin{prop}
Soit $X$ un plongement projectif de $G/H$. Un diviseur anticanonique de $X$ est $$-K_X=\sum_{i=1}^m X_i+\sum_{\alpha\in S\backslash I}a_\alpha D_\alpha$$ o\`u $a_\alpha$ est l'entier $\langle2\rho^P,\check{\alpha}\rangle$ et $2\rho^P$ est le caract\`ere donn\'e par $\sum_{\alpha\in R^+\backslash R_I}\alpha$.
\end{prop}
Pour prouver ce r\'esultat (qui est un corollaire de \cite[th.4.2]{Br97a}), on se ram\`ene au cas o\`u $X$ est un plongement toro\"idal. On note $\pi:X\longrightarrow G/P$ la fibration de fibre torique $Y$. On a alors $K_X=\pi^*(K_{G/P})+K_\pi$. La fin de la preuve utilise le r\'esultat connu dans le cas torique $-K_X=\sum_{i=1}^m X_i$ et le r\'esultat pour les vari\'et\'es de drapeaux $$-K_{G/P}=\sum_{\alpha\in S\backslash I}a_\alpha D_\alpha.$$

Par construction de l'\'eventail colori\'e associ\'e \`a $X$, les diviseurs irr\'educti\-bles $G$-stables $X_i$ correspondent aux ar\^etes de l'\'eventail $\mathbb F$ associ\'e \`a $X$ qui ne sont pas engendr\'ees par une couleur. On note $x_i\in N$ l'\'el\'ement primitif de l'ar\^ete correspondante \`a $X_i$.\\

Afin de d\'eterminer les plongements de Fano de $G/H$, on utilise la caract\'erisation suivante des diviseurs amples \cite[th.3.3]{Br89}.
\begin{prop}\label{cardivample}
Soient $X$ un plongement projectif de $G/H$, $\mathbb F$ son \'eventail colori\'e et $D$ un diviseur de Weil de la forme $$\sum_{i=1}^m b_i X_i+\sum_{\alpha\in S\backslash I}b_\alpha D_\alpha,$$ o\`u les $b_i$ et les $b_\alpha$ sont dans $\Zbb$.

Alors $D$ est de Cartier si et seulement si pour tout c\^one colori\'e $(\mathcal{C},\mathcal{F})$ de $\mathbb F$, il existe $\chi_{\mathcal{C}}$ dans $M$ tel que $$\forall x_i \in\mathcal{C},\,\langle x_i,\chi_C\rangle=b_i$$  $$\mbox{et }\forall D_\alpha\in\mathcal{F},\,\langle\check\alpha_M,\chi_{\mathcal{C}}\rangle =b_\alpha.$$

Lorsque $D$ est de Cartier on peut alors d\'efinir une application $h_D$ de $N_\Rbb$ dans $\Rbb$, comme suit. Soit $x$ un \'el\'ement de $N_\Rbb$. Comme $X$ est projectif, $\mathbb F$ est complet donc il existe un unique c\^one colori\'e maximal $(\mathcal{C},\mathcal{F})$ de $\mathbb F$ tel que $x$ soit dans $\mathcal{C}$, et on pose alors $h_D(x)=\langle x,\chi_{\mathcal{C}}\rangle$. C'est une application lin\'eaire sur chaque c\^one.

Le diviseur (de Cartier) $D$ est ample si et seulement si:\\
(i) l'application $h_D$ est strictement convexe (c'est-\`a-dire, pour chaque c\^one maximal $\mathcal{C}$, l'application lin\'eaire $\chi_{\mathcal{C}}$ est strictement sup\'erieure sur $\mathcal{C}$ aux applications lin\'eaires  $\chi_{\mathcal{C'}}$, pour tout c\^one maximal $\mathcal{C'}$ distinct de $\mathcal{C}$);\\
(ii) pour tout c\^one colori\'e $(\mathcal{C},\mathcal{F})$ de $\mathbb F$ et pour tout $\alpha$ tel que $D_\alpha$ ne soit pas dans $\mathcal{F}$, on a $$\langle\check\alpha_M,\chi_{\mathcal{C}}\rangle <b_\alpha.$$
\end{prop}

Petits rappels (pour plus de d\'etails, voir la partie \ref{notations}): $S$ est l'ensemble des racines simples et $I$ est le sous-ensemble de $S$ de la d\'efinition \ref{INM}. Les r\'eseaux $N$ et $M$ sont des r\'eseaux de rang \'egal au rang de $G/H$, et sont duaux l'un de l'autre. Pour tout $\alpha\in S\backslash I$, $\check\alpha_M$ est l'image de la couleur $D_\alpha$ dans $N$ par l'application $\sigma$ (\ref{sigma}).

On d\'efinit $\mathcal{D}_X$ comme l'ensemble des $\alpha\in S\backslash I$ tels que $D_\alpha$ soit une couleur de $X$.\\

Lorsque $D=-K_X$ est ample, l'application $h_D$ v\'erifie $h_D(x_1)=\cdots=h_D(x_m)=1$ pour tout $i\in\{1,\dots,m\}$, et $h_D(\check\alpha_M)=a_\alpha$ pour tout $\alpha\in\mathcal{D}_X$. L'ensemble
\begin{equation}\label{defipolytope}
Q=\{u\in N_\Rbb\mid \forall \mathcal{C}\in\mathbb{F},\,\langle\chi_{\mathcal{C}},u\rangle\leq 1\}
\end{equation}
est alors un polytope convexe contenant $0$ dans son int\'erieur et dont les sommets sont $x_1,\dots,x_m$ et certains des $\frac{\check\alpha_M}{a_\alpha}$ o\`u $\alpha\in\mathcal{D}_X$. De plus les points $\frac{\check\alpha_M}{a_\alpha}$ o\`u $\alpha\not\in\mathcal{D}_X$ sont dans l'int\'erieur de $Q$, et les $\frac{\check\alpha_M}{a_\alpha}$ o\`u $\alpha\in\mathcal{D}_X$ sont dans le bord de $Q$. En particulier, $Q$ est un polytope convexe rationnel.

On d\'efinit le polytope dual de $Q$ par
\begin{equation}
Q^*=\{v\in M_\Rbb\mid \forall u\in Q,\,\langle v,u\rangle\geq-1\}.
\end{equation}
Les sommets de $Q^*$ sont en bijection avec les faces maximales de $Q$. Ainsi $Q^*$ est l'enveloppe convexe des $-\chi_{\mathcal{C}}$ o\`u $\mathcal{C}$ d\'ecrit l'ensemble des c\^ones colori\'es maximaux de $\mathbb{F}$. En particulier, $Q^*$ est \`a sommets dans $M$.

Je regroupe toutes ces propri\'et\'es dans la d\'efinition suivante.
\begin{defi} \label{reflexif}
Soit $G/H$ un espace  homog\`ene horosph\'erique.
Un polytope convexe $Q$ de $N_\Rbb$ est dit  {\it $G/H$-r\'eflexif} si les trois conditions suivantes sont v\'erifi\'ees:\\

(1) $Q$ est \`a sommets dans $N\cup\{\frac{\check\alpha_M}{a_\alpha}\mid\alpha\in S\backslash I\}$, et contient 0 dans son int\'erieur.\\

(2) $Q^*$ est \`a sommets dans $M$.\\

(3) Pour tout $\alpha\in S\backslash I$, $\frac{\check\alpha_M}{a_\alpha}\in Q$.\\
\end{defi}

\begin{exs}
D\'esormais on repr\'esentera, dans les figures, les couleurs par les $\frac{\check\alpha_M}{a_\alpha}$ au lieu des $\check\alpha_M$, et toujours par des points blancs. Les exemples choisis seront tels qu'un point blanc est associ\'e \`a une seule couleur. Mais il faut tout de m\^eme noter qu'il peut y avoir des cas o\`u un m\^eme point $\frac{\check\alpha_M}{a_\alpha}$ correspond \`a plusieurs couleurs.\\

(1) Dans le cas de $SL_2/U$, on n'a qu'une couleur et $a_\alpha=2$. Les seuls polytopes $SL_2/U$-r\'eflexifs sont les segments suivants.
\begin{center}
\includegraphics[width=13cm]{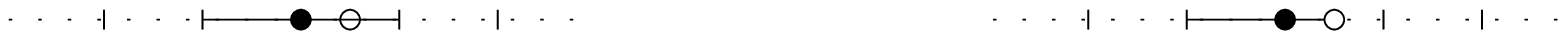}
\end{center}
Et leurs duaux sont respectivement:
\begin{center}
\includegraphics[width=13cm]{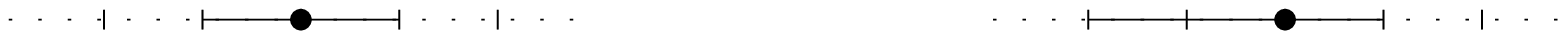}
\end{center}

(2)\label{polyref} Dans le cas de $(SL_2\times\Cbb)/U$, prenons les polytopes associ\'es aux \'eventails de l'exemple \ref{SL2C}. Les deux premiers sont $G/H$-r\'eflexifs. Le troisi\`eme polytope n'est pas $G/H$-r\'eflexif car son dual n'est pas \`a sommets entiers.

\begin{center}
\includegraphics[width=13cm]{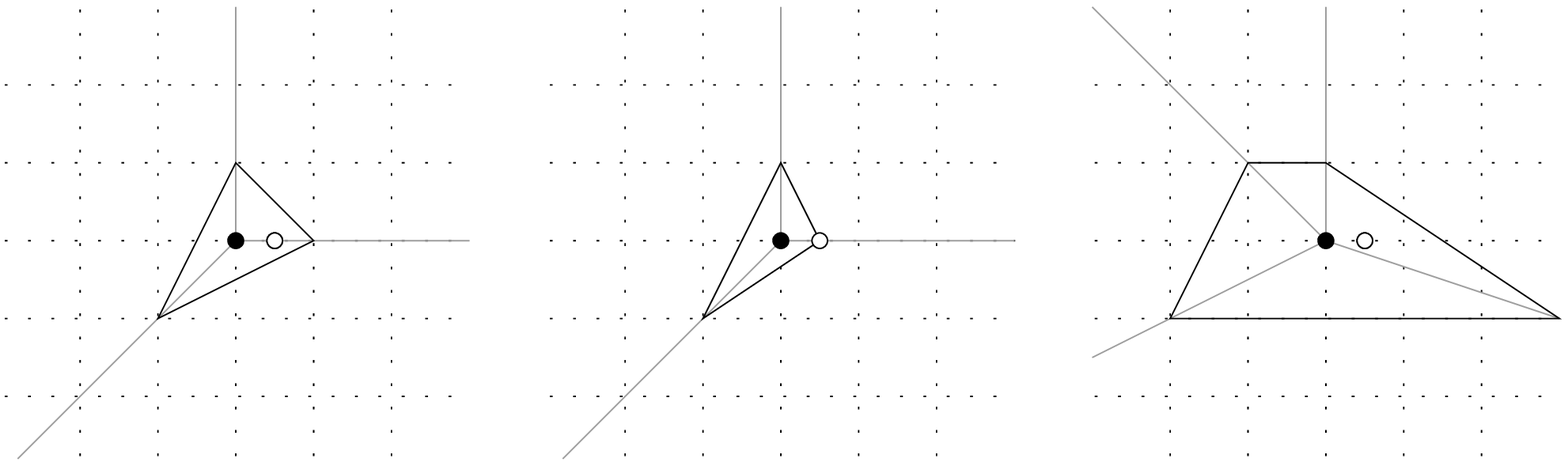}
\end{center}
Les polytopes duaux sont respectivement:
\begin{center}
\includegraphics[width=13cm]{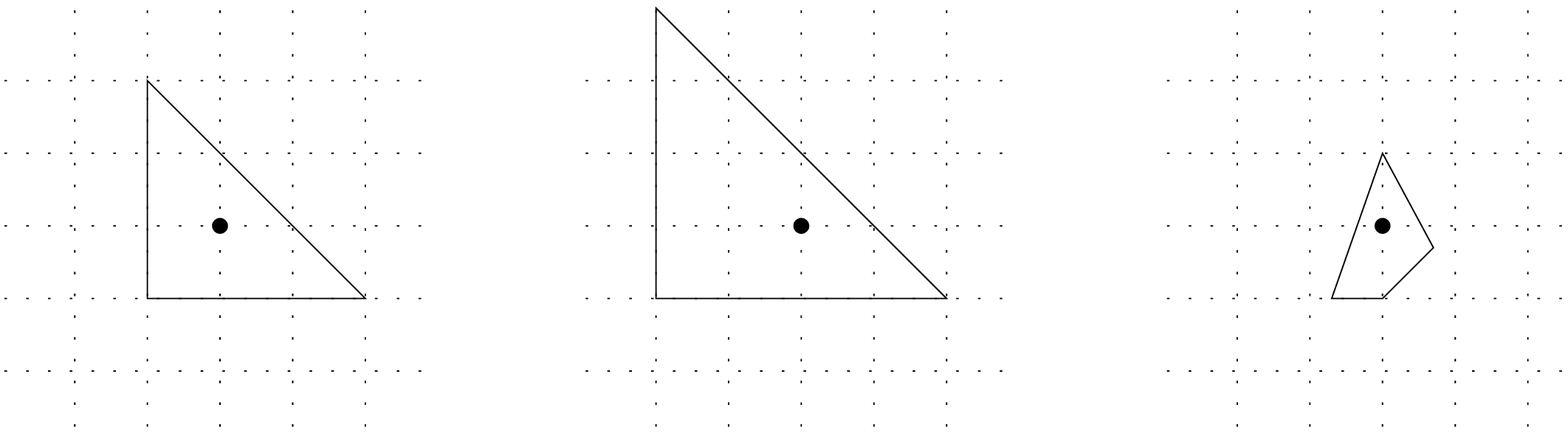}
\end{center}

(3) Voici d'autres exemples de polytopes dans le cas o\`u $G/H=(SL_2\times\Cbb)/U$.
\begin{center}
\includegraphics[width=13cm]{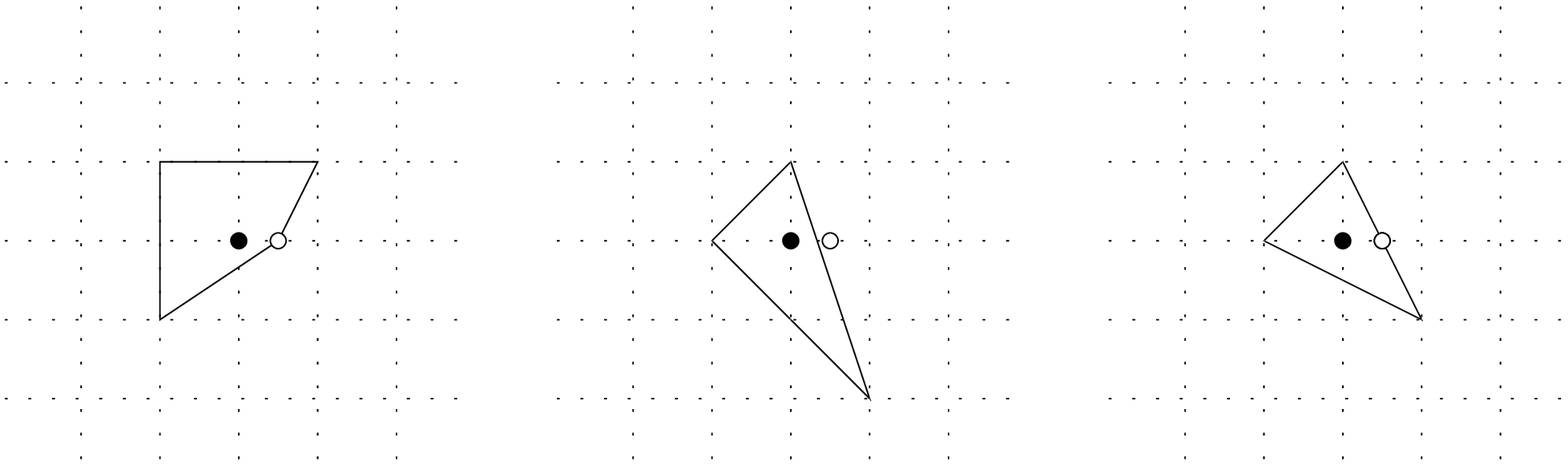}
\end{center}
Le premier et le troisi\`eme sont $G/H$-r\'eflexifs. Par contre, le deuxi\`eme ne l'est pas, car $\frac{\check\alpha_M}{a_\alpha}\not\in Q$.
\end{exs}

\begin{rem}\label{polymoment}
La condition (3) est \'equivalente \`a la condition \\
(3'): $2\rho^P+Q^*$ (le translat\'e  de $Q^*$ par $2\rho^P$ dans $\Lambda_\Rbb$) est inclus dans $\Lambda^+_\Rbb$.

En fait, $2\rho^P+Q^*$ est le polytope moment de $(X,-K_X)$ \cite{Br89}, c'est-\`a-dire
$$H^0(X,-K_X)=\bigoplus_{\lambda\in Q^*}V(2\rho^P+\lambda).$$
En effet, si $s$ est la section canonique de $-K_X$, on a l'isomorphisme suivant
$$\begin{array}{ccc}
\{f\in\Cbb(G/H)^{(B)}\mid \operatorname{div}(f)-K_X\geq 0\} & \lra & H^0(X,-K_X)^{(B)} \\
  f & \lmt & fs.
\end{array}$$
De plus $s$ est de poids $2\rho^P=\sum_{\alpha\in S\backslash I}a_\alpha\omega_\alpha$ sous l'action de $B$, et si $f\in\Cbb(G/H)^{(B)}$ est de poids $\chi$ alors $$\operatorname{div}(f)-K_X=\sum_{i=1}^m (1+\langle\chi,x_i\rangle)X_i +\sum_{\alpha\in S\backslash I}(a_\alpha+\langle\chi,\check\alpha_M\rangle)D_\alpha.$$
Donc les poids de $H^0(X,-K_X)^{(B)}$ sont bien les caract\`eres qui sont dans le polytope $2\rho^P+Q^*$.
\end{rem}

Il peut \^etre int\'eressant aussi de regarder les cas o\`u $-K_X$ est seulement $\Qbb$-Cartier (c'est-\`a-dire un multiple de $-K_X$ est de Cartier) et ample: on dira alors que $X$ est $\Qbb$-Fano. Dans ce cas, $Q$ ne v\'erifie plus la condition (2). On dit alors qu'un polytope convexe de $N_\Rbb$ est $\Qbb$-$G/H$-r\'eflexif s'il v\'erifie la condition (3) de la d\'efinition \ref{reflexif} et la condition\\
(1'): les sommets de $Q$ sont des \'el\'ements primitifs de $N$ ou des \'el\'ements de $\{\frac{\check\alpha_M}{a_\alpha},\alpha\in S\backslash I\}$, et $Q$ contient 0 dans son int\'erieur.

On peut remarquer qu'un polytope $G/H$-r\'eflexif v\'erifie aussi la condition (1)' car (1) et (2) impliquent (1').

Le dernier polytope de l'exemple \ref{polyref} (2) est  $\Qbb$-$G/H$-r\'eflexif.\\

Remarquons que dans le cas torique, la d\'efinition d'un polytope $(\Cbb^*)^n$-r\'eflexif est celle d'un polytope r\'eflexif donn\'ee par V.~Batyrev.
L'ensemble des classes d'isomorphisme des vari\'et\'es toriques de Fano de dimension $n$ est en bijection avec l'ensemble des polytopes r\'eflexifs de $\Zbb^n$ \cite{Ba94}. La proposition suivante g\'en\'eralise alors cette classification aux vari\'et\'es horosph\'eriques de Fano.

\begin{prop}
Soit $G/H$ un espace homog\`ene horosph\'erique.
L'application qui \`a un plongement $X$ de Fano de $G/H$ associe le polytope $Q$ d\'efini en (\ref{defipolytope}) est une bijection de l'ensemble des plongements de Fano de $G/H$ (respectivement $\Qbb$-Fano) \`a isomorphisme pr\`es, sur l'ensemble des polytopes $G/H$-r\'eflexifs (respectivement $\Qbb$-$G/H$-r\'eflexifs) de $N_\Rbb$.
\end{prop}
Je rappelle que les isomorphismes de plongements sont d\'efinis en \ref{plongiso}.
\begin{proof}
Il suffit de d\'efinir l'application inverse. Soit $Q$ un polytope $G/H$-r\'eflexif de $N_\Rbb$; on lui associe alors le plongement $X(Q)$ de $G/H$ dont l'\'eventail colori\'e est l'ensemble des c\^ones colori\'es $(\mathcal{C},\mathcal{F})$ (et leurs faces colori\'ees) tels que $\mathcal{C}$ soit le c\^one engendr\'e par une face maximale $F$ de $Q$, et $\mathcal{F}$ soit l'ensemble des $D_\alpha$ v\'erifiant $\frac{\check\alpha_M}{a_\alpha}\in F$.

Alors on v\'erifie que $-K_{X(Q)}$ est de Cartier ($-\chi_\mathcal{C}$ est le sommet de $Q^*$ associ\'e \`a $F$) et ample, par  convexit\'e de $Q$ et par la condition (3) de la d\'efinition \ref{reflexif}. Ensuite, le polytope associ\'e \`a $X(Q)$ (\ref{defipolytope}) est bien $Q$ par construction.

De m\^eme, si $X$ est un plongement de Fano et $Q$ est le polytope $G/H$-r\'eflexif associ\'e, on a $X(Q)=X$.

Lorsque $Q$ est un polytope $\Qbb$-$G/H$-r\'eflexif, on construit $X(Q)$ exactement de la m\^eme fa\c{c}on, et on v\'erifie aussi que cette application est l'inverse de l'application $X\mapsto Q(X)$.
\end{proof}
\begin{rem}\label{couleurbord}
Les couleurs d'un plongement $X$ de Fano sont les couleurs $D_\alpha$ telles que $\frac{\check\alpha_M}{a_\alpha}$ se trouve sur le bord du polytope $G/H$-r\'eflexif associ\'e \`a $X$.
Supposons que pour $\alpha\neq\beta$ on ait $\frac{\check\alpha_M}{a_\alpha}=\frac{\check\beta_M}{a_\beta}$. Soit $X$ un plongement de Fano de $G/H$. Alors si $\frac{\check\alpha_M}{a_\alpha}$ est sur le bord du polytope $G/H$-r\'eflexif associ\'e \`a $X$ alors $\frac{\check\beta_M}{a_\beta}$ l'est aussi; et r\'eciproquement. Par cons\'equent, soit aucune des deux couleurs $D_\alpha$ et $D_\beta$ n'est une couleur de $X$, soit toutes les deux sont des couleurs de $X$.
\end{rem}
Cette classification permet notamment de donner une version effective, dans le cas des vari\'et\'es horosph\'eriques de Fano, d'un r\'esultat de V. Alexeev et M. Brion \cite{AB04} sur les vari\'et\'es sph\'eriques de Fano.
\begin{teo}\label{finitude}
Soit $G/H$ un espace homog\`ene horosph\'erique de rang n. On note $a=\prod_{\alpha\in S\backslash I}a_\alpha$, et $V=(7(a+1))^{n2^{n+1}}$.
L'ensemble des classes d'isomorphisme des vari\'et\'es  de Fano qui sont des plongements de $G/H$ est fini  et de cardinal inf\'erieur \`a $$(n!aV)^{\frac{n(n+1)}{2}}2^{2^n(n!aV)^{n+1}}.$$
\end{teo}
\begin{rem}
La borne obtenue est sans doute loin d'\^etre optimale. En effet, les majorations faites dans la preuve sont en g\'en\'eral assez grossi\`eres.
\end{rem}

On appelle alors {\it automorphisme de $(N,\mathcal{D})$} tout automorphisme du r\'eseau $N$ qui fixe chaque couleur $\check\alpha_M$.
Gr\^ace \`a la proposition suivante, il suffira de majorer le nombre de polytopes $G/H$-r\'eflexifs \`a automorphisme de $(N,\mathcal{D})$ pr\`es pour d\'emontrer le th\'eor\`eme  \ref{finitude}.

\begin{prop}\label{autom}
Soit $\phi$ un automorphisme de $(N,\mathcal{D})$. On note aussi $\phi$ l'automorphisme de $N_\Rbb$ induit par $\phi$. Soient $X$ et $X'$ des plongements de $G/H$ d'\'eventails colori\'es respectifs $\mathbb{F}$ et $\mathbb{F'}$, tels que $\mathbb{F'}=\phi(\mathbb{F})$. Alors les vari\'et\'es $X$ et $X'$ sont isomorphes.
\end{prop}
Rappelons que les plongements $X$ et $X'$ ne sont pas isomorphes si $\phi $ n'est pas trivial (th\'eor\`eme \ref{LunaVust}).

La preuve de cette proposition est fortement inspir\'ee de \cite{AB04}.
\begin{proof}
Notons $G'$ la partie semi-simple de $G$. On d\'efinit alors $\tilde{G}=G'\times P/H$. C'est un groupe alg\'ebrique r\'eductif et connexe.  Rappelons que $P=N_G(H)$, ainsi $\tilde{G}$ agit sur $G/H$ par $(g,pH).xH=gxpH$.

De plus $G/H$ est homog\`ene sous l'action de $\tilde{G}$, $G/H\simeq \tilde{G}/\tilde{H}$, et $\tilde{H}$ est un sous-groupe horosph\'erique de $\tilde{G}$. En effet, le radical unipotent $\tilde{U}$ du sous-groupe de Borel $\tilde{B}=(B\cap G')\times P/H$ est $U\times\{1\}$ et $$\tilde{H}=\{(g,pH)\in\tilde{G}\mid gp\in H\}\supset \tilde{U}.$$
On v\'erifie aussi que $\tilde{P}:=N_{\tilde{G}}(\tilde{H})=(P\cap G')\times P/H$. Puis on voit que le r\'eseau des caract\`eres de $\tilde{P}$ dont la restriction \`a $\tilde{H}$ est triviale est le groupe $\{(\chi,\chi^{-1}),\,\chi\in M\}$ isomorphe \`a $M$. De plus, les couleurs associ\'ees \`a $\tilde{G}/\tilde{H}$ sont les m\^emes que celles associ\'ees \`a $G/H$ avec les m\^emes images dans $M$.

Ainsi les plongements de $G/H$ sont les m\^emes que ceux de $\tilde{G}/\tilde{H}$. Donc $X$ et $X'$ sont des plongements de $\tilde{G}/\tilde{H}$ d'\'eventails colori\'es respectifs $\mathbb{F}$ et $\mathbb{F'}$.

Comme $P/H$ est isomorphe au tore dual de $M$, $\phi$ induit naturellement un automorphisme de $P/H$, qu'on notera encore $\phi$. On peut alors d\'efinir $X''$ comme \'etant la vari\'et\'e $X$ sur laquelle $\tilde{G}$ agit par $(g,pH)_\phi .x=gx\phi^{-1}(pH)$; c'est un plongement de $\tilde{G}/\tilde{H}$. Il suffit alors de montrer que son \'eventail colori\'e est $\mathbb{F'}$: on aura ainsi des isomorphismes de vari\'et\'es $X'\simeq X''\simeq X$.

Notons $\phi^*:M\lra M$ l'automorphisme dual de $\phi$. Comme $\phi$ fixe chaque couleur, $\phi^*(\chi)$ et $\chi$ ont alors la m\^eme restriction \`a $B\cap G'$, pour tout $\chi\in M$.

Soit $f\in\Cbb(\tilde{G}/\tilde{H})^{(\tilde{B})}$ de poids $\chi$; autrement dit, pour tous $(b,pH)\in\tilde{B}$ et $x\in\tilde{G}/\tilde{H}$, on a $f(bxpH)=\chi^{-1}(b)\chi (p)f(x)$.

Si l'action de $\tilde{G}$ est tordue par $\phi$, on a alors $$f((b,pH)_\phi .x)=\chi^{-1}(b)\chi (\phi^{-1}(pH))f(x).$$ Or $\chi (\phi^{-1}(pH))=\phi^{*-1}(\chi)(p)$ par d\'efinition de $\phi^*$ et $\chi^{-1}(b)=\phi^{*-1}(\chi^{-1})(b)$ car $\phi$ fixe chaque couleur. Donc $f$ est de poids $\phi^{*-1}(\chi)$ (lorsque l'action de $\tilde{G}$ est tordue par $\phi$).

Par cons\'equent, pour tout diviseur $\tilde{B}$-stable $D$ de $X$, en notant $D''$ le m\^eme diviseur dans $X''$, on a $\sigma(D'')=\phi(\sigma(D))$. De la fa\c{c}on dont l'\'eventail colori\'e d'un plongement est construit \`a partir des images par $\sigma$ des diviseurs $B$-stables, on conclut que $X''$ a pour \'eventail colori\'e $\mathbb{F'}$.
\end{proof}

\begin{ex}
On sait d\'ej\`a que le nombre de vari\'et\'es toriques de Fano de dimension $2$ est $16$, dont $5$ seulement sont lisses.

En rang $2$, il y a par exemple $135$ polytopes $(SL_2\times\Cbb^*)/U$-r\'eflexifs, dont $16$ correspondent \`a un plongement lisse et $398$ polytopes $(SL_2\times SL_2)/U$-r\'eflexifs, dont $39$ correspondent \`a un plongement lisse. Les polytopes $SL_3/U$-r\'eflexifs sont les m\^emes que les polytopes $(SL_2\times SL_2)/U$-r\'eflexifs; par contre, seulement $27$ d'entre eux correspondent \`a un plongement lisse de Fano de $SL_3/U$. La liste compl\`ete de tous ces polytopes est donn\'ee dans \cite[ch.6]{Pa06}.
\end{ex}

Dans le cas torique, le th\'eor\`eme \ref{finitude} dit que le nombre de polytopes r\'eflexifs, \`a automorphisme de $\Zbb^n$ pr\`es, est fini. Ce r\'esultat a \'et\'e d\'emontr\'e en premier par A.~Borisov et L.~Borisov \cite{BB92}. Mais une autre preuve (plus souvent utilis\'ee) consiste \`a appliquer un r\'esultat de D.~Hensley \cite{He83} afin de majorer le volume des polytopes r\'eflexifs.

Dans notre cas, on utilise une g\'en\'eralisation du th\'eor\`eme de D.~Hensley due \`a J.~Lagarias et G.~Ziegler. On note $\stackrel{o}{Q}$ l'int\'erieur de $Q$.
\begin{teo} [\cite{LZ91}]\label{LZ}
Soient $n$, $a$ et $k$ des entiers strictement positifs. Soit $Q\subset \Rbb^n$ un polytope convexe \`a sommets dans $\Zbb^n$ tel que $\sharp (\stackrel{o}{Q}\cap a\Zbb^n)=k$. Alors le volume de $Q$ est major\'e par $ka^n(7(ka+1))^{n2^{n+1}}$.
\end{teo}

\begin{proof}[D\'emonstration du th\'eor\`eme \ref{finitude}]
Commen\c{c}ons par remarquer qu'un polytope $G/H$-r\'eflexif $Q$ est \`a sommets dans $\frac{1}{a}N$ et que $\stackrel{o}{Q}\cap N=\{0\}$. En effet, s'il existe $u$ dans $(\stackrel{o}{Q}\cap N)\backslash\{0\}$, alors pour tout $v\in Q^*$, $\langle v,u\rangle>-1$. Or il existe un sommet $v$ de $Q^*$ tel que $\langle v,u\rangle<0$, ce qui contredit le fait que $\langle v,u\rangle$ est entier.

Par cons\'equent, le th\'eor\`eme pr\'ecedent nous dit que tout polytope $G/H$-r\'eflexif a un volume major\'e par $V=(7(a+1))^{n2^{n+1}}$.

Si $G$ est semi-simple, les $\frac{\check\alpha_M}{a_\alpha}$ pour $\alpha\in S\backslash I$ forment une famille g\'en\'eratrice de $N_\Rbb$.
Soient $\alpha_1,\dots,\alpha_n\in S\backslash I$ tels que $(\check{\alpha}_{1M},\dots,\check{\alpha}_{nM})$ soit une base de $N_\Rbb$. Soit $u\in Q$, alors pour tout $i\in\{1,\dots,n\}$, le simplexe de sommets $0$, $u$, et $\frac{\check{\alpha}_{jM}}{\alpha_j},\,j\neq i$ est strictement inclus  dans $Q$, donc son volume est strictement inf\'erieur \`a $V$. Par cons\'equent, $u$ est dans l'int\'erieur du parall\'el\'epip\`ede de sommets $\pm n!aV\check{\alpha}_{1M},\dots,\pm n!aV\check{\alpha}_{nM}$. Ce dernier contient au maximum $2^n(n!aV)^{n+1}$ points de $N$, donc en majorant grossi\`erement, on obtient que le nombre de polytopes $G/H$-r\'eflexifs est, dans ce cas, inf\'erieur \`a $2^{2^n(n!aV)^{n+1}}$.

Dans le cas g\'en\'eral o\`u $G$ est r\'eductif, notons $N_\Rbb^1$ le sous-espace vectoriel de $N_\Rbb$ engendr\'e par les $\check\alpha_M$. Soit $(\check\alpha_{1M},\dots,\check\alpha_{lM})$ une base de $N_\Rbb^1$, alors pour tout polytope $G/H$-r\'eflexif $Q$, il existe $f_1,\dots,f_{n-l}\in Q\cap N$ tels que $\{\check\alpha_{1M},\dots,\check\alpha_{lM})\}\cup\{f_1,\dots,f_{n-l}\}$ forme une base de $N_\Rbb$. Montrons qu'il existe un ensemble fini fix\'e d'\'el\'ements de $N_\Rbb$, tel qu'on puisse toujours se ramener, quitte \`a appliquer un automorphisme de $(N,\mathcal{D})$, au cas o\`u  les $f_i$ font partie de cet ensemble.\\
Soit $(e_1,\dots,e_n)$ une base fix\'ee de $N$ telle que $(e_1,\dots,e_l)$ soit une base de $N^1=N_\Rbb^1\cap N$. Dans cette base, un automorphisme de $(N,\mathcal{D})$ s'\'ecrit alors matriciellement sous la forme $\mathcal A=\left(
\begin{array}{rr}
I_l & B\\
0 & C
\end{array}
\right)$, o\`u $I_l$ est l'identit\'e dans $N^1$, $B\in M_{l,n-l}(\Zbb)$ et $C\in GL_{n-l}(\Zbb)$. Soit $\mathcal F$ la matrice form\'ee des vecteurs colonnes $f_1,\dots,f_{n-l}$ dans la base $(e_1,\dots,e_n)$. Alors, il existe un automorphisme $\mathcal A$ de $(N,\mathcal{D})$ tel que $\mathcal G=\mathcal A \mathcal F$ de la forme suivante (on note  $(g_{ij})_{1\leq i\leq n,1\leq j\leq n-l}$ ses coefficients)
$$ \mathcal G=\left(
\begin{array}{rrrr}
g_{11} & &  \cdots & g_{1,n-l}\\
\vdots & & & \vdots\\
g_{l+1,1}& & & \vdots\\
0 & \ddots & & \vdots\\
\vdots & \ddots & \ddots & \vdots\\
0 & \cdots & 0 & g_{n,n-l}
\end{array}
\right) $$
et telle que pour tout $1\leq j\leq n-l$ et $1\leq i\leq j-1$, $$0\leq g_{ij}<|g_{j+l,j}|.$$
De plus, le produit des $g_{j+l,j}$ est major\'e par $n!aV$, donc le nombre de telles matrices $\mathcal G$ est inf\'erieur \`a $(n!aV)^{\frac{n(n+1)}{2}}$.

On en d\'eduit alors que le nombre de polytopes $G/H$-r\'eflexifs est inf\'erieur \`a $$(n!aV)^{\frac{n(n+1)}{2}}2^{2^n(n!aV)^{n+1}}.$$
\end{proof}
\begin{rem}
Ce r\'esultat peut se g\'en\'eraliser de la mani\`ere suivante: soit $\kappa$ un entier strictement positif, alors le nombre de plongements $X$ de $G/H$, tels que $-\kappa K_X$ soit de Cartier et ample, est fini. Il suffit de remarquer que le polytope $\Qbb$-$G/H$-r\'eflexif $Q$ associ\'e \`a $X$ a son dual \`a sommets dans $\frac{1}{\kappa}M$. On en d\'eduit alors que $\frac{1}{\kappa}Q$ est \`a sommets dans $\frac{1}{a\kappa}N$, et son int\'erieur ne contient que $0$ comme point entier. Ainsi le volume de $Q$ est major\'e par $\kappa^n(7(a\kappa+1))^{n2^{n+1}}$, et la suite de la d\'emonstration reste inchang\'ee.
\end{rem}

\section{Majoration du degr\'e et du nombre de Picard}

Pour majorer le degr\'e et le nombre de Picard des vari\'et\'es horosph\'eriques de Fano, on se donne un espace homog\`ene horosph\'erique $G/H$ et on \'etudie le degr\'e et le nombre de Picard de ses plongements lisses.

\subsection{Majoration du degr\'e}
Pour avoir plus de d\'etails sur le degr\'e des vari\'et\'es de Fano, on se r\'ef\'erera \`a \cite{De03} o\`u le degr\'e est major\'e dans le cas des vari\'et\'es toriques lisses, et aussi  \`a \cite{De01} o\`u on trouve en particulier des vari\'et\'es toriques de Fano de grand degr\'e (prop.5.22).

\begin{defi}
Soit $X$ une vari\'et\'e de Fano de dimension $d$. On appelle {\it degr\'e} de $X$, le nombre d'intersection $(-K_X)^d$. Le th\'eor\`eme de Riemann-Roch et le th\'eor\`eme d'annulation de Serre impliquent que la dimension de $\Gamma(X,-kK_X)$ est \'equivalente, lorsque l'entier $k$ tend vers l'infini, \`a $\frac{k^d}{d!\,}(-K_X)^d$.
\end{defi}
Le r\'esultat suivant se d\'emontre en utilisant la remarque \ref{polymoment} et la formule des caract\`eres de Weyl.
\begin{prop}\label{degre}(cas particulier de \cite[th.4.1]{Br89})
Soit $X$ un plongement de Fano de $G/H$. Alors $$(-K_X)^d=d!\,\int_{Q^*}\prod_{\alpha\in R^+\backslash R_I^+}\frac{\langle 2\rho^P+\chi,\check\alpha\rangle}{\langle\rho^B,\check\alpha\rangle}d\chi,$$
o\`u la mesure dans l'int\'egrale est celle pour laquelle le domaine fondamental de $M$ est de volume~$1$.
\end{prop}

Dans le cas torique, on a $(-K_X)^d=d!\,\operatorname{vol}(Q^*)$ \cite[cor 2.23]{Od88}. Le th\'eor\`eme  \ref{LZ} permet alors de majorer le degr\'e de $X$ de fa\c{c}on imm\'ediate, mais cette borne est doublement exponentielle. O.~Debarre donne une bien meilleure borne dans le cas o\`u $X$ est lisse.
\begin{teo}[\cite{De03}]\label{degretorique}
Soit $X$ une vari\'et\'e torique lisse de Fano de dimension $n$ et de nombre de Picard $\rho$.\\
Si $\rho>1$, alors $(-K_X)^n\leq n!\,n^{\rho n}$.
\end{teo}
Si $\rho=1$, on a bien s\^ur $X=\Pbb^n$ et $(-K_X)^n=(n+1)^n$.
Dans le cas horosph\'erique, on obtient le r\'esultat analogue \'enonc\'e dans l'introduction: le th\'eor\`eme \ref{majorationdegre}.

La proposition \ref{locfac} se r\'e\'ecrit de la fa\c{c}on suivante lorsque $X$ est de Fano.
\begin{prop} \label{locfac2} Soit $G/H$ un espace homog\`ene horosph\'erique. Soit $X$ un plongement de Fano de $G/H$.  On note $Q$ le polytope $G/H$-r\'eflexif associ\'e. Rappelons que $\mathcal{D}_X$ d\'esigne l'ensemble des \'el\'ements de $S\backslash I$ qui correspondent \`a une couleur de $X$. Alors $X$ est localement factoriel si et seulement si pour toute face maximale $F$ de $Q$ on a:\\
(1) les seuls points de $F$ de la forme $\frac{\check\alpha_M}{a_\alpha}$ sont des sommets de $F$;\\
(2) $F$ est un simplexe dont les sommets $e_1,\dots,e_n$ v\'erifient:\\
\indent (i) pour tout $i$, ou bien $e_i\in N$ ou bien il existe un unique $\alpha\in\mathcal{D}_X$ tel que $e_i=\frac{\check\alpha_M}{a_\alpha}$,\\
\indent (ii) $(a_1e_1,\dots,a_ne_n)$ est une base de $N$ o\`u les $a_i$ sont d\'efinis par:
$$\begin{array}{lcl}
a_i=1 & \mbox{si} & e_i\in N\\
a_i=a_\alpha & \mbox{si} & e_i=\frac{\check\alpha_M}{a_\alpha}.
\end{array}$$
\end{prop}
\begin{rem}
Lorsque $Q$ est le polytope $G/H$-r\'eflexif associ\'e \`a un plongement de Fano localement factoriel de $G/H$, on d\'efinira $a_u$ comme dans la proposition pour tout sommet $u$ de~$Q$.
\end{rem}
Dans la suite, $G/H$ sera un espace homog\`ene horosph\'erique fix\'e et $X$ un plongement de Fano localement factoriel de $G/H$.

On note $r$ l'entier strictement positif tel que le nombre de sommets de $Q$ soit \'egal \`a $n+r$.

Dans le cas torique, cet entier est le nombre de Picard. Dans le cas horosph\'erique, on peut exprimer le nombre de Picard de la fa\c{c}on suivante:
\begin{equation}\label{equationrho}
\rho=m+\sharp (S\backslash I)-n=r+\sharp (S\backslash I)-\sharp(\mathcal{D}_X).
\end{equation}

En effet, tout diviseur est lin\'eairement \'equivalent \`a un diviseur de la forme $\sum_{i=1}^m b_i X_i+\sum_{\alpha\in S\backslash I}b_\alpha D_\alpha$ ($b_i$ et $b_\alpha$ entiers) \cite[Chap 5]{Br97b}, et les relations entre ces diviseurs sont les relations du type $\operatorname{div} (f)=0$ o\`u $f\in\Cbb(G/H)^{(B)}$. Pour la deuxi\`eme \'egalit\'e, rappelons que $\sharp(\mathcal{D}_X)$ est le nombre de couleurs de $X$, c'est-\`a-dire le nombre des points $\frac{\check\alpha_M}{a_\alpha}$ qui sont sur le bord de $Q$ (voir la remarque \ref{couleurbord}). Puisque $X$ est localement factoriel, les points de la forme $\frac{\check\alpha_M}{a_\alpha}$ sont tous des sommets de $Q$ et correspondent \`a une unique couleur. On en d\'eduit que $\sharp(\mathcal{D}_X)$ est exactement le nombre de sommets de la forme $\frac{\check\alpha_M}{a_\alpha}$. Ainsi $m+\sharp(\mathcal{D}_X)=n+r$.

De plus, ces deux \'egalit\'es restent vraies lorsque $X$ est seulement $\Qbb$-factoriel (voir la d\'efinition \ref{qlocfac}). Elles se d\'emontrent par les m\^emes arguments et avec la proposition \ref{carqlocfac}.

\begin{ex}\label{r=1}
Si $r=1$, alors $Q$ est un simplexe: c'est l'enveloppe convexe de $e_1,\dots,e_{n+1}$ et pour tout $i\in\{1,\dots,n+1\}$, \\ $(a_1e_1,\dots, a_{i-1}e_{i-1},a_{i+1}e_{i+1},\dots,a_{n+1}e_{n+1})$ est une base de $N$. Donc $$a_{n+1}e_{n+1}=-a_1e_1-\cdots-a_ne_n.$$

Par la suite on verra que le cas o\`u $r=1$ est un cas un peu \`a part, comme dans le cas torique.
\end{ex}

\begin{prop}\label{vol}
Soient $X$ un plongement de Fano localement factoriel de $G/H$ et $Q$ le polytope $G/H$-r\'eflexif associ\'e.
On note $$C=n+\sum_{\alpha\in S\backslash I}(a_\alpha -1).$$
Si $\rho\geq 2$, alors $$\operatorname{vol}(Q^*)\leq(C^r\max_{\alpha\in S\backslash I}a_\alpha)^n,$$
et si $\rho=1$, on a $$\operatorname{vol}(Q^*)\leq((C+1)\max_{\alpha\in S\backslash I}a_\alpha)^n.$$
\end{prop}
On peut observer que $C$ ne d\'epend que de $G$ et $H$. On verra par la suite (lemme \ref{lem3}) que $C\leq d$.

La preuve de cette proposition est inspir\'ee de la preuve du th\'eor\`eme \ref{degretorique} \cite{De03}. Elle consiste en deux lemmes.
\begin{lem}\label{lem0}
Soit $b$ un r\'eel strictement positif. Si pour tout  $u\in Q$ et pour tout $v\in Q^*$, on a $-1\leq \langle v,u\rangle\leq b$,
alors  $$\operatorname{vol}(Q^*)\leq((b+1)\max_{\alpha\in S\backslash I}a_\alpha)^n.$$
\end{lem}
\begin{proof}
Par d\'efinition du dual d'un polytope, on a $-\frac{1}{b}Q^*\subset Q^*$: en effet si $v\in -\frac{1}{b}Q^*$, on a bien $\langle v,u\rangle\geq -1$ pour tout $u\in Q$ par hypoth\`ese.

Raisonnons par l'absurde et supposons que $\operatorname{vol}(Q^*)>((b+1)\max_{\alpha\in S\backslash I}a_\alpha)^n$.\\
Soient un r\'eel $\eta\in]0,1[$, et $$Q'=\frac{1-\eta}{(b+1)\max_{\alpha\in S\backslash I}a_\alpha}Q^*.$$ Par hypoth\`ese, le volume de $Q'$ est strictement sup\'erieur \`a $1$ pour $\eta$ assez petit. Donc, par le th\'eor\`eme de van der Corput \cite{Co36}, il existe deux \'el\'ements $q$ et $q'$ de $Q$ tels que $q-q'$ soit dans $M$. De plus $Q'$ est convexe et $-\frac{1}{b}Q'\subset Q'$, donc $$\frac{q-q'}{b+1}=\frac{q-b(-\frac{1}{b}q')}{b+1}\in Q',\\ \mbox{ et } q-q'\in (b+1)Q'\subset\frac{1}{\max_{\alpha\in S\backslash I}a_\alpha}\stackrel{o}{Q^*}.$$
Il reste \`a montrer que $0$ est le seul point de $M$ dans l'int\'erieur de $\frac{1}{\max_{\alpha\in S\backslash I}a_\alpha}Q^*$ pour obtenir la contradiction voulue.

Soit $v\in M$ non nul dans l'int\'erieur de $\frac{1}{\max_{\alpha\in S\backslash I}a_\alpha}Q^*$. Pour tout $u\in Q$, on a $$\langle v,u\rangle>-\frac{1}{\max_{\alpha\in S\backslash I}a_\alpha}.$$ Or il existe un sommet $u$ de $Q$ tel que $\langle v,u\rangle$ soit strictement n\'egatif, et ce sommet est dans $\frac{1}{a_\alpha}\Zbb$ pour un $\alpha\in S\backslash I$ donn\'e. Ceci n'est pas possible, donc un tel $v$ n'existe pas.
\end{proof}
\begin{ex}
Revenons \`a l'exemple \ref{r=1} o\`u $r=1$; on a pour tout $i\in\{1,\dots,n+1\}$, $a_i(1+\langle v_j,e_i\rangle)=\delta_{ij}(a_1+\cdots+a_{n+1})$ o\`u $v_j$ est le sommet de $Q^*$ associ\'e \`a la face de $Q$ oppos\'ee \`a $e_j$, et $\delta_{ij}$ est le symbole de Kronecker. Dans ce cas, on peut prendre $b=a_1+\cdots+a_{n+1}-1$. Regardons ce que vaut $b$ selon les valeurs de $\rho$. Sachant que $\rho$ vaut $1$ plus le nombre de racines simples $\alpha\in S\backslash I$ telles que $\frac{\check\alpha_M}{a_\alpha}$ ne soit pas un sommet de $Q$,  quand $\rho>1$, la somme $a_1+\cdots+a_{n+1}$ peut \^etre major\'ee par $C$. Par contre, si $\rho=1$,  $a_1+\cdots+a_{n+1}$ vaut exactement $C+1$. On obtient dans les deux cas le r\'esultat de la proposition \ref{vol} lorsque $r=1$.
\end{ex}
Le lemme suivant g\'en\'eralise le r\'esultat obtenu dans l'exemple ci-dessus.
\begin{lem}\label{lem1}
Si $r\geq 2$, alors pour tout sommet $u$ de $Q$ et tout sommet $v$ de $Q^*$, on a $$\,0\leq a_u(1+\langle v,u\rangle)\leq C^r.$$
\end{lem}

\begin{proof}
Soit $v$ un sommet de $Q^*$. On note $e_1,\dots,e_n$ les sommets de la face de $Q$ associ\'ee \`a $v$. Alors $(a_1e_1,\dots,a_ne_n)$ est une base de $N$, et $\langle v,e_i\rangle=-1$ pour tout $i\in\{1,\dots,n\}$.

On note $0<b_1<\cdots<b_k$, $k\leq r$, les \'el\'ements de l'ensemble suivant: $$\{a_u(1+\langle v,u\rangle)\mid u\mbox{ sommet de }Q\mbox{ distinct des }e_i\}.$$ On va alors montrer par r\'ecurrence sur $j$ que $b_j$ est major\'e par $C^j$.

Soient $j\in\{1,\dots,k\}$ et $u$ un sommet de $Q$ tel que $a_u(1+\langle v,u\rangle)=b_j$. Soit $E$ un sous-ensemble de $\{1,\dots,n\}$ non vide, minimal, tel que $u$ et les $(e_i)_{i\in E}$ ne soient pas sur une face commune de $Q$. Un tel sous-ensemble existe car $u$ et  $e_1,\dots,e_n$ ne sont pas sur la m\^eme face. Quitte \`a changer l'ordre des $e_i$, supposons $E=\{1,\dots,\epsilon\}$. Posons $w=a_uu+a_1e_1+\cdots+a_\epsilon e_\epsilon$, alors $w$ est un \'el\'ement de $N$.

Il existe une face de $Q$ telle que $w$ soit dans le c\^one engendr\'e par cette face. De plus, cette face n'est pas celle correspondant \`a $v$. Il existe donc $u_1,\dots,u_s$, parmi les sommets de cette face, distincts de $e_1,\dots,e_n$ (mais pas forc\'ement distincts entre eux), et $e'_1,\dots,e'_t$ parmi $e_1,\dots,e_n$ (toujours pas forc\'ement distincts entre eux), tels que
\begin{equation}
w=\sum_{i=1}^sa_{u_i}u_i+\sum_{i=1}^ta'_ie'_i.\label{relationetoile}
\end{equation}
Posons $a'_i=a_{e'_i}$ pour simplifier les notations.

La relation (\ref{relationetoile}) n'est pas triviale, car $u$ et les $(e_i)_{i\in E}$ ne sont pas sur une face commune de $Q$. De plus, pour tout $i\in\{1,\dots,t\}$, $e'_i$ n'est pas dans $\{e_1,\dots,e_\epsilon\}$. En effet, si par exemple $e'_t=e_\epsilon$, alors la relation (\ref{relationetoile}) induit une relation la relation non triviale suivante: $$a_uu+a_1e_1+\cdots+a_{\epsilon-1}e_{\epsilon-1}=\sum_{i=1}^sa_{u_i}u_i+\sum_{i=1}^{t-1}a'_ie'_i.$$ On en d\'eduit alors que $u, e_1,\dots,e_{\epsilon-1}$ ne sont pas sur une face commune de $Q$, ce qui contredit la minimalit\'e de $E$.

On va utiliser la relation \ref{relationetoile} pour majorer $b_j$ en fonction de $b_{j-1}$. Pour cela, on va d'abord \'etablir une relation entre les entiers $a_u$, $a_{u_i}$, $a_i$ et $a'_i$.

Le point $$\frac{1}{a_u+\sum_{i=1}^\epsilon a_i}w$$ est le barycentre des points $u,e_1,\dots,e_\epsilon$ affect\'es des coefficients (strictement positifs) respectifs $a_u,a_1,\dots,a_\epsilon$. Comme $u$ et  $e_1,\dots,e_\epsilon$ sont des sommets de $Q$ qui ne font pas partie d'une m\^eme face, ce barycentre est dans l'int\'erieur de $Q$. Par cons\'equent, $$\frac{\sum_{i=1}^sa_{u_i}+\sum_{i=1}^ta'_i}{a_u+\sum_{i=1}^\epsilon a_i}<1,$$ autrement dit, $$a_u+\sum_{i=1}^\epsilon a_i\geq\sum_{i=1}^sa_{u_i}+\sum_{i=1}^ta'_i+1.$$
La relation \ref{relationetoile} nous donne alors $$b_j=a_u(1+\langle v,u\rangle)=\sum_{i=1}^sa_{u_i}(1+\langle v,u_i\rangle)+a_u+\sum_{i=1}^\epsilon a_i-\sum_{i=1}^sa_{u_i}-\sum_{i=1}^ta'_i.$$
Si $s=0$ alors $b_j=a_u+\sum_{i=1}^\epsilon a_i-\sum_{i=1}^ta'_i\leq a_u+\sum_{i=1}^\epsilon a_i$.

Sinon, $b_j\geq\sum_{i=1}^sa_{u_i}(1+\langle v,u_i\rangle)+1$ et donc $a_{u_i}(1+\langle v,u_i\rangle)<b_j$ pour tout $i\in\{1,\dots,s\}$; autrement dit, $a_{u_i}(1+\langle v,u_i\rangle)\leq b_{j-1}$. Ainsi,
\begin{eqnarray*}
b_j & \leq & sb_{j-1}+a_u+\sum_{i=1}^\epsilon a_i-\sum_{i=1}^sa_{u_i}-\sum_{i=1}^ta'_i \\
 & \leq & sb_{j-1}+a_u+\sum_{i=1}^\epsilon a_i-s.
 \end{eqnarray*}
Or $s\leq\sum_{i=1}^sa_{u_i}\leq\sum_{i=1}^sa_{u_i}+\sum_{i=1}^ta'_i\leq a_u+\sum_{i=1}^\epsilon a_i$, donc $$b_j\leq s(b_{j-1}-1)+a_u+\sum_{i=1}^\epsilon a_i\leq (a_u+\sum_{i=1}^\epsilon a_i)b_{j-1}.$$
En particulier si $j=1$, alors $s=0$ et donc $b_1\leq a_u+\sum_{i=1}^\epsilon a_i$.

Pour conclure, il suffit juste de remarquer que si $\epsilon=n$ on a alors $r=1$. Donc $\epsilon\leq n-1$ et \\$$a_u+\sum_{i=1}^\epsilon a_i\leq n+\sum_{\alpha\in S\backslash I}(a_\alpha -1)=C$$ car les entiers $a_u$ et $a_i$ valent soit $1$ soit $a_\alpha$.
\end{proof}
La proposition se d\'eduit alors facilement de ces deux lemmes: le lemme \ref{lem1} donne le r\'eel $b$ \`a utiliser dans le lemme \ref{lem0}.

Pour montrer le th\'eor\`eme \ref{majorationdegre}, il suffit maintenant d'une part de majorer le terme \`a l'int\'erieur de l'int\'egrale de la proposition \ref{degre}, et d'autre part de donner une borne explicite pour $C$.
\begin{lem}\label{lem2}
Si $r\geq 2$, alors pour tout $\alpha\in R^+\backslash R_I^+$ et tout $\chi\in Q^*$, $$0\leq\frac{\langle 2\rho^P+\chi,\check\alpha\rangle}{\langle\rho^B,\check\alpha\rangle}\leq C^r\max_{\alpha\in S\backslash I}a_\alpha.$$
\end{lem}
\begin{proof}
Il suffit de montrer le r\'esultat pour tout sommet $\chi$ de $Q^*$.

Soit $\alpha\in S\backslash I$.

Si $\frac{\check\alpha_M}{a_\alpha}$ est un sommet de $Q$, alors le lemme \ref{lem1} nous dit que $$\frac{\langle 2\rho^P+\chi,\check\alpha\rangle}{\langle\rho^B,\check\alpha\rangle}=a_\alpha(1+\langle\chi,\frac{\check\alpha_M}{a_\alpha}\rangle)\leq C^r.$$

Si $\frac{\check\alpha_M}{a_\alpha}$ n'est pas un sommet de $Q$, on obtient une majoration l\'eg\`erement plus grande. Comme $\frac{\check\alpha_M}{a_\alpha}\in Q$, il existe des sommets $u_i$ de $Q$ et des r\'eels positifs $\lambda_i$ tels que $\sum\lambda_i\leq 1$ et $\frac{\check\alpha_M}{a_\alpha}=\sum\lambda_i u_i$. D'apr\`es le lemme \ref{lem1} on a pour tout $i$, $\langle\chi,u_i\rangle\leq C^r-1$ car $a_{u_i}\geq 1$. On en d\'eduit alors que $$1+\langle\chi,\frac{\check\alpha_M}{a_\alpha}\rangle\leq C^r\mbox{ et donc }\frac{\langle 2\rho^P+\chi,\check\alpha\rangle}{\langle\rho^B,\check\alpha\rangle}\leq C^ra_\alpha.$$

Soit maintenant $\alpha\in R^+\backslash R_I^+$. Ecrivons $\alpha=\beta_1+\cdots+\beta_s$ o\`u $\beta_1,\dots,\beta_s$ sont des racines simples. On a alors $$\frac{\langle 2\rho^P+\chi,\check\alpha\rangle}{\langle\rho^B,\check\alpha\rangle}=\frac{\sum_{i=1}^s\langle 2\rho^P+\chi,\check{\beta_i}\rangle}{\sum_{i=1}^s1}\leq\frac{sC^r\max_{\alpha\in S\backslash I}a_\alpha}{s}.$$
\end{proof}
\begin{rem}
Lorsque $\rho=1$, d'une part $r=1$ et d'autre part $\frac{\check\alpha_M}{a_\alpha}$ est un sommet de $Q$ pour tout $\alpha\in S\backslash I$. On a alors le m\^eme r\'esultat en rempla\c{c}ant $C^r\max_{\alpha\in S\backslash I}a_\alpha$ par $C+1$.
\end{rem}
\begin{lem}\label{lem3}
 $$\sum_{\alpha\in S\backslash I}(a_\alpha -1)\leq\dim(G/P)=\sharp(R^+\backslash R_I^+).$$
En particulier $C\leq d=n+\sharp(R^+\backslash R_I^+)$.
\end{lem}
\begin{proof}
La preuve se fait en \'etudiant les diff\'erents cas, mais avant tout, r\'eduisons le nombre de ces cas.\\

Etape 1. Pour tout $\alpha\in S\backslash I$, $a_\alpha=\langle 2\rho^P,\check\alpha\rangle=\langle 2\rho^B,\check\alpha\rangle-\langle 2\rho_I,\check\alpha\rangle$ o\`u $2\rho_I=\sum_{\beta\in R_I^+}\beta$. Ainsi, comme $\langle 2\rho^B,\check\alpha\rangle=2$, on a $$\sum_{\alpha\in S\backslash I}(a_\alpha -1)=\sharp(S\backslash I)-\sum_{\alpha\in S\backslash I}\sum_{\beta\in R_I^+}\langle\beta,\check\alpha\rangle.$$
Le r\'esultat \`a montrer est donc $$\sum_{\alpha\in S\backslash I}\sum_{\beta\in R_I^+}-\langle\beta,\check\alpha\rangle\leq\sharp(R^+\backslash R_I^+)-\sharp(S\backslash I).$$

Etape 2. Montrons qu'il suffit de montrer le r\'esultat lorsque le diagramme de Dynkin $\Gamma_I$ (voir la d\'efinition \ref{lisse}) est connexe. Supposons donc le r\'esultat vrai dans ce cas. Soit $I=\bigsqcup_{j=1}^tI_j$ tel que $\Gamma_I=\bigsqcup_{j=1}^t\Gamma_{I_j}$ soit la d\'ecomposition de $\Gamma_I$ en composantes connexes. Notons $S_j$ l'ensemble des racines simples de $S$ li\'ees \`a $I_j$ dans le diagramme de Dynkin $\Gamma_S$ de $G$ (c'est-\`a-dire les racines simples $\alpha$ telles qu'il existe une racine simple $\beta$ de $I_j$ avec $\langle\alpha,\check\beta\rangle\not=0$). Alors $$\sum_{\alpha\in S\backslash I}\sum_{\beta\in R_I^+}-\langle\beta,\check\alpha\rangle=\sum_{j=1}^t\sum_{\beta\in R_{I_j}^+}\sum_{\alpha\in S_j\backslash I_j}-\langle\beta,\check\alpha\rangle$$ ce qui, par hypoth\`ese, est inf\'erieur \`a $\sum_{j=1}^t(\sharp(R_{S_j}^+\backslash R_{I_j}^+)-\sharp(S_j\backslash I_j))$.\\

Or $\sharp(R_{S_j}^+\backslash R_{I_j}^+)-\sharp(S_j\backslash I_j)$ est le nombre de racines positives non simples de $R_{S_j}^+\backslash R_{I_j}^+$. De plus, $R_{S_j}^+\cap R_{S_k}^+=S_j\cap S_k$ si $j\not=k$ car deux racines simples ne peuvent pas \^etre toutes les deux li\'ees \`a $I_j$ et $I_k$ (un diagramme de Dynkin ne contient pas de cycle). On a donc $$\sum_{j=1}^t(\sharp(R_{S_j}^+\backslash R_{I_j}^+)-\sharp(S_j\backslash I_j))\leq\sharp(R^+\backslash R_I^+)-\sharp(S\backslash I).$$

Il suffit donc de montrer le r\'esultat lorsque $\Gamma_I$ est connexe et $S$ est l'ensemble des racines simples li\'ees \`a $I$.\\

Etape 3. On peut de plus supposer que $S\backslash I$ est un singleton. En effet,  en regardant la r\'eunion disjointe $\bigsqcup_{\alpha\in S\backslash I}R_{I\cup\{\alpha\}}^+\backslash R_I^+$ incluse dans $R^+\backslash R_I^+$, on voit que $$\sum_{\alpha\in S\backslash I}(\sharp(R_{I\cup\{\alpha\}}^+\backslash R_I^+)-1)\leq\sharp(R^+\backslash R_I^+)-\sharp(S\backslash I).$$

Etape 4. Il reste \`a calculer d'une part $-\sum_{\beta\in R_I^+}\langle\beta,\check\alpha\rangle$, et d'autre part $\sharp(R_{I\cup\{\alpha\}}^+\backslash R_I^+)-1$ dans les diff\'erents cas o\`u $\Gamma_I$ est connexe et $S\backslash I$ est un singleton.

Le tableau ci-dessous r\'esume le r\'esultat des calculs en fonction du type des diagrammes de Dynkin $\Gamma_I$ et $\Gamma_S$. On peut remarquer qu'il y a deux fa\c{c}ons diff\'erentes d'ajouter un sommet \`a un diagramme de type $A_1$ pour obtenir un diagramme de type $G_2$. Il y a aussi deux fa\c{c}ons d'ajouter un sommet \`a un diagramme de type $A_1$ pour obtenir un  diagramme de type $B_2$;  l'une des fa\c{c}on est compt\'ee dans le cas $A_i$ et $B_{i+1}$ et l'autre dans le cas $A_i$ et $C_{i+1}$ en identifiant $B_2$ et $C_2$.

\begin{center}
\begin{tabular}{||l||c|c||}
\hline
Types respectifs de $\Gamma_I$ et de $\Gamma_S$ & $ -\sum_{\beta\in R_I^+}\langle\beta,\check\alpha\rangle$ & $\sharp(R_{I\cup\{\alpha\}}^+\backslash R_I^+)-1$\\
\hline
\hline
$A_i$ et $A_{i+1}$ pour $i\geq 0$ & $i$ & $i$\\
\hline
$A_i$ et $B_{i+1}$ pour $i\geq 1$ & $2i$ & $\frac{i(i+3)}{2}$\\
\hline
$A_i$ et $C_{i+1}$ pour $i\geq 1$ & $i$ & $\frac{i(i+3)}{2}$\\
\hline
$A_i$ et $D_{i+1}$ pour $i\geq 3$ & $2(i-1)$ & $\frac{(i-1)(i+2)}{2}$\\
\hline
$A_5$ et $E_6$ & $9$ & $20$\\
\hline
$A_6$ et $E_7$ & $12$ & $41$\\
\hline
$A_7$ et $E_8$ & $15$ & $91$\\
\hline
$A_1$ et $G_2$ (fl\`eche vers $\alpha$) & $3$ & $4$\\
\hline
$A_1$ et $G_2$ (fl\`eche partant de $\alpha$) & $1$ & $4$\\
\hline
\hline
$B_i$ et $B_{i+1}$ pour $i\geq 2$ & $2i-1$ & $2i$\\
\hline
$B_2$ et $C_3$ & $4$ & $4$\\
\hline
$B_3$ et $F_4$ & $9$ & $14$\\
\hline
\hline
$C_i$ et $A_{i+1}$ pour $i\geq 3$ & $2i$ & $2i$\\
\hline
$C_3$ et $F_4$ & $6$ & $14$\\
\hline
\hline
$D_i$ et $D_{i+1}$ pour $i\geq 4$ & $2i-2$ & $2i-1$\\
\hline
$D_5$ et $E_6$ & $10$ & $15$\\
\hline
$D_6$ et $E_7$ & $15$ & $32$\\
\hline
$D_7$ et $E_8$ & $21$ & $78$\\
\hline
\hline
$E_6$ et $E_7$ & $16$ & $26$\\
\hline
$E_7$ et $E_8$ & $27$ & $56$\\
\hline
\end{tabular}
\end{center}

\begin{exs}

1/ Lorsque $\Gamma_I$ est de type $A_i$ et que $\Gamma_S$ est de type $A_{i+1}$, notons $\alpha,\beta_1,\dots,\beta_i$ les racines simples de $S$,
\begin{center}\includegraphics{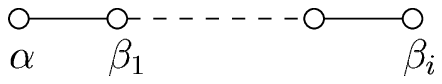}\end{center}  alors $$-\sum_{\beta\in R_I^+}\langle\beta,\check\alpha\rangle= -\sum_{1\leq j\leq k\leq i}\langle\beta_j+\cdots+\beta_k,\check\alpha\rangle$$ $$=-\sum_{1\leq k\leq i}\langle\beta_1+\cdots+\beta_k,\check\alpha\rangle=-\sum_{1\leq k\leq i}-1=i.$$ Puis $R_{I\cup\{\alpha\}}^+\backslash R_I^+=\{\alpha+\beta_1+\cdots +\beta_k\mid 1\leq k\leq i\}$.

2/ Lorsque $\Gamma_I$ est de type $A_1$ et que $\Gamma_S$ est de type $G_2$ (fl\`eche vers $\alpha$), notons $\beta_1$ la racine de~$I$
\begin{center}\includegraphics{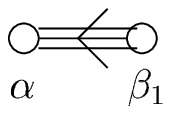}\end{center} alors $$-\sum_{\beta\in R_I^+}\langle\beta,\check\alpha\rangle=-\langle\beta_1,\check\alpha\rangle=3.$$
Puis $R_{I\cup\{\alpha\}}^+\backslash R_I^+=\{\alpha,\alpha +\beta,2\alpha +\beta,3\alpha +\beta,3\alpha +2\beta\}$.
\end{exs}
\end{proof}
On peut remarquer que les configurations o\`u il y a \'egalit\'e sont celles qu'on retrouve lorsqu'on s'int\'eresse aux vari\'et\'es horosph\'eriques non toro\"idales lisses (voir la partie \ref{sectionlisse}). Le r\'esultat de ce lemme est alors optimal m\^eme dans le cas lisse.\\

On a maintenant tous les outils pour d\'emontrer le th\'eor\`eme \ref{majorationdegre}.
\begin{proof}
Si $\rho>1$ alors les propositions \ref{degre} et \ref{vol} ainsi que les lemmes \ref{lem2} et \ref{lem3} nous permettent de dire que $$(-K_X)^d\leq d!\,(d^r\max_{\alpha\in S\backslash I}a_\alpha)^n(d^r\max_{\alpha\in S\backslash I}a_\alpha)^{d-n}=d!\,d^{rd}(\max_{\alpha\in S\backslash I}a_\alpha)^d.$$

On a \'evidemment $\max_{\alpha\in S\backslash I}a_\alpha\leq C\leq d$ d\`es que $n\geq 1$. D'autre part, si $n=0$ on a $X=G/P$ et  $\rho=\sharp(S\backslash I)>1$; donc on a aussi $\max_{\alpha\in S\backslash I}a_\alpha\leq C$ puisque $a_\alpha\geq 2$ pour tout $\alpha\in S\backslash I$.

Si $\rho\geq r+1$ alors $(-K_X)^d\leq d!\,d^{d\rho}$.

Si $\rho=r>1$, $\frac{\check\alpha_M}{a_\alpha}$ est un sommet de $Q$ pour tout $\alpha\in S\backslash I$. La preuve du lemme \ref{lem2} nous dit que $$\frac{\langle 2\rho^P+\chi,\check\alpha\rangle}{\langle\rho^B,\check\alpha\rangle}\leq C^r$$ et donc $$(-K_X)^d\leq d!\,(d^r\max_{\alpha\in S\backslash I}a_\alpha)^n(d^r)^{d-n}=d!\,d^{rd}(\max_{\alpha\in S\backslash I}a_\alpha)^n\leq d!\,d^{d\rho+n}.$$

De m\^eme, si $\rho=1$ on a aussi $r=1$ et $$(-K_X)^d\leq d!\,((d+1)\max_{\alpha\in S\backslash I}a_\alpha)^n((d+1))^{d-n}=d!\,(d+1)^{d}(\max_{\alpha\in S\backslash I}a_\alpha)^n\leq d!\,(d+1)^{d+n}.$$
\end{proof}

\subsection{Majoration du nombre de Picard}
On va majorer le nombre de Picard des vari\'et\'es horosph\'eriques $\Qbb$-factori\-elles.
\begin{defi}\label{qlocfac}
Une vari\'et\'e normale est dite {\it $\Qbb$-factorielle} si tout diviseur de Weil est $\Qbb$-Cartier.
\end{defi}
On va formuler un crit\`ere analogue \`a celui donn\'e dans la proposition \ref{locfac2}, toujours en utilisant la caract\'erisation des diviseurs de Cartier sur une vari\'et\'e sph\'erique \cite[prop.3.1]{Br89} \'enonc\'ee dans la proposition \ref{cardivample}. La preuve est laiss\'ee au lecteur.
\begin{prop}\label{carqlocfac}
Soit $X$ un plongement de Fano de $G/H$. On note $Q$ le polytope $G/H$-r\'eflexif associ\'e. Alors $X$ est $\Qbb$-factoriel si et seulement si toute face $F$ de $Q$ est un simplexe, et tous les points de $F$ de la forme $\frac{\check\alpha_M}{a_\alpha}$ sont des sommets de $F$.
\end{prop}
On va montrer un r\'esultat analogue \`a celui obtenu par C.~Casagrande dans le cas torique \cite[th.1(i)]{Ca06}.
\begin{teo}\label{picardnumber}
Soit $X$ une vari\'et\'e horosph\'erique de Fano, $\Qbb$-factorielle, de rang $n$, de dimension $d$ et de nombre de Picard $\rho$. Alors
$$\rho\leq 2n+\sharp(S\backslash I)\leq n+d\leq 2d.$$
\end{teo}
On en d\'eduit alors facilement le th\'eor\`eme \ref{nombrepicard} \'enonc\'e dans l'introduction, en remarquant que $n+d=2d$ si et seulement si $X$ est torique.

La d\'emonstration qui suit est inspir\'ee des preuves de C.~Casagrande \cite[th.3(i)]{Ca06} et de B.~Nill \cite[lem.5.5]{Ni05}.
\begin{proof}
Soit $v$ un sommet de $Q^*$.

Etape 1. Notons $F_v$ la face de $Q$ associ\'ee \`a $v$, et $e_1,\dots,e_n$ les sommets de cette face. Montrons que tout sommet entier $u$ (c'est-\`a-dire dans $N$) de $Q$ v\'erifiant $\langle v,u\rangle=0$ est adjacent \`a $F_v$, c'est-\`a-dire qu'il existe un indice $j$ tel que $e_1,\dots,e_{j-1},e_{j+1},\dots,e_n$ et $u$ soient les sommets d'une face $F_j$ de $Q$. On note $u^j$ le sommet $u$ de $Q$ v\'erifiant cette derni\`ere condition.

Soit $(e_1^*,\dots,e_n^*)$ la base duale de $(e_1,\dots,e_n)$ dans $M_\Rbb$. Soit $j\in\{1,\dots,n\}$. Alors $\langle e_j^*,u^j\rangle\neq 0$ (sinon $u^j$ est dans l'hyperplan engendr\'e par les $(e_i)_{i\neq j}$) et on peut alors d\'efinir $$\gamma_j=\frac{-1-\langle v,u^j\rangle}{\langle e_j^*,u^j\rangle}.$$
De plus, le sommet de $Q^*$ associ\'e \`a $F_j$ est $v^j=v+\gamma_je_j^*$. On en d\'eduit alors que $\gamma_j>0$, car $\langle v^j,e_j\rangle >-1$.

Soit $u$ un sommet entier de $Q$ v\'erifiant $\langle v,u\rangle=0$; alors $\langle v^j,u\rangle=\gamma_j\langle e^*_j,u\rangle$ et donc $$u\not\in F_j \Longleftrightarrow \langle e_j^*,u\rangle\geq 0.$$

Si  $u\neq u^j$ pour tout $j\in\{1,\dots,n\}$, alors  $\langle e_j^*,u\rangle\geq 0$ pour tout $j\in\{1,\dots,n\}$ et donc $u$ est dans le c\^one engendr\'e par $e_1,\dots,e_n$, ce qui n'est pas possible. Par cons\'equent, $u$ est l'un des $u^j$ et est donc adjacent \`a $F_v$.\\

Etape 2. Le nombre de sommets de $Q$ tels que $\langle v,u\rangle=-1$ est $n$, et le nombre de sommets entiers de $Q$ tels que $\langle v,u\rangle=0$ est inf\'erieur ou \'egal \`a $n$ par l'\'etape 1.

L'origine est dans $Q$, donc il existe des sommets $v_1,\dots,v_h$ de $Q^*$ ($h>0$) et des entiers strictement positifs $m_1,\dots,m_h$ tels que $m_1v_1+\cdots +v_hm_h=0$.

On note $I=\{1,\dots,h\}$, $M=\sum_{i\in I}m_i$, et pour tout sommet entier $u$ de $Q$ on pose
\begin{eqnarray*}
  A(u)&=&\{i\in I\mid\langle v_i,u\rangle=-1\}\\
 \mbox{ et }  B(u)&=&\{i\in I\mid\langle v_i,u\rangle=0\}.
 \end{eqnarray*}
Alors on a, pour tout sommet entier $u$ de $Q$,
\begin{eqnarray*}
0 & = & \sum_{i\in I}m_i\langle v_i,u\rangle\\
  & = & -\sum_{i\in A(u)}m_i+\sum_{i\not\in A(u)\cup B(u)}m_i\langle v_i,u\rangle\\
  & \geq & -\sum_{i\in A(u)}m_i+\sum_{i\not\in A(u)\cup B(u)}m_i\\
  & = & M-2\sum_{i\in A(u)}m_i-\sum_{i\in B(u)}m_i
\end{eqnarray*}
et donc $M\leq 2\sum_{i\in A(u)}m_i+\sum_{i\in B(u)}m_i$.

Sommons cette derni\`ere in\'egalit\'e sur tous les sommmets entiers $u$ de $Q$. On obtient alors, en notant $r'$ le nombre de ces sommets:
\begin{eqnarray*}
r'M & \leq &\sum_u\sum_{i\in A(u)}2m_i+\sum_u\sum_{i\in B(u)}m_i\\
    & = & \sum_{i\in I}\sum_{u,\langle v_i,u\rangle=-1}2m_i+\sum_{i\in I}\sum_{u,\langle v_i,u\rangle=0}m_i\\
    & \leq  & 3nM
\end{eqnarray*}
donc le nombre de sommets entiers de $Q$ est inf\'erieur ou \'egal \`a $3n$. On en d\'eduit alors facilement que le nombre $n+r$ de sommets  de $Q$ est inf\'erieur ou \'egal \`a $3n+\sharp(\mathcal{D}_X)$. En utilisant l'\'equation (\ref{equationrho}), on a $$\rho=r+\sharp (S\backslash I)-\sharp(\mathcal{D}_X)\leq 2n+\sharp (S\backslash I).$$
\end{proof}
\begin{cor}
Soit $X$ une vari\'et\'e horosph\'erique de Fano localement factorielle de dimension $d\geq 2$. Alors  $$(-K_X)^d\leq d!\,d^{3d^2}.$$
\end{cor}

\begin{proof}
Lorsque $\rho\geq 2$ le degr\'e est major\'e par $$d!\,d^{d(2n+\sharp(S\backslash I))+n}\leq d!\,d^{d(2n+\sharp(S\backslash I)+n)},$$ gr\^ace aux th\'eor\`emes \ref{majorationdegre} et \ref{picardnumber}. Sachant que $n+\sharp(S\backslash I)\leq d$ (\ref{dimension}), on obtient facilement le r\'esultat.

Lorsque $\rho=1$, on a $(-K_X)^d\leq d!\,(d+1)^{d+n}\leq d!\,(d+1)^{2d}\leq d!\,d^{3d^2}$ d\`es que $d>1$.
\end{proof}
\begin{rems}
Le cas o\`u $d=1$ et $\rho=1$ correspond \`a la droite projective; le degr\'e vaut alors~2.

Si on fait le m\^eme raisonnement dans le cas torique (en combinant le th\'eor\`eme de Debarre et celui de Casagrande), on obtient une borne l\'eg\`erement plus forte, $d!\,d^{2d^2}$ au lieu de $d!\,d^{3d^2}$.
\end{rems}

\section{Sur l'amplitude des diviseurs d'une vari\'et\'e horosph\'erique projective}\label{chapample}
G. Ewald et U. Wessels ont montr\'e que pour toute vari\'et\'e torique projective $X$ de dimension~$d$, et pour tout diviseur de Cartier ample $D$, le diviseur $(d-1)D$ est tr\`es ample \cite{EW91}. On va utiliser leur argument pour d\'emontrer le th\'eor\`eme \ref{tresample}.

On peut supposer que $G$ est factoriel, quitte \`a le remplacer par le produit direct de son radical et du rev\^etement de sa partie semi-simple. Alors tout fibr\'e en droites $\mathcal{L}$ sur une $G$-vari\'et\'e normale $X$ est $G$-lin\'earisable (voir \cite[d\'ef.2.1 et ch.2.4]{KK89}), et donc $H^0(X,\mathcal{L})$ est un $G$-module rationnel par \cite[lem.2.5]{KK89}.

Soit $X$ une $G$-vari\'et\'e sph\'erique projective de rang $n$, et soit $D$ un diviseur de Cartier ample sur $X$. On peut supposer que $D$ est $B$-stable, puisque tout diviseur sur une vari\'et\'e sph\'erique est lin\'eairement \'equivalent \`a un diviseur $B$-stable.

Soit $Q_D^*$ le polytope moment de $D$; c'est un polytope convexe d'int\'erieur non vide dans $M_\Rbb$, d\'efini par $$ Q_D^*=\{\lambda\in\Lambda_\Qbb\mid\mbox{ le }G\mbox{-module }V(k\lambda)\mbox{ appara\^it dans }H^0(X,kD)\mbox{ pour un certain entier }k>0\}.$$ En particulier, $$H^0(X,D)\simeq \bigoplus_{\lambda\in Q_D^*}V(\lambda).$$

Lorsque $X$ est horosph\'erique et $D=-K_X$, on retrouve le polytope moment $2\rho^P+Q^*$ de la remarque \ref{polymoment}.

Toujours lorsque $X$ est horosph\'erique, le polytope $Q_D^*$ est enti\`erement d\'efini par l'espace des sections de $D$. En effet, notons $s$ la section canonique de $D$. C'est un vecteur propre sous l'action de $B$. On note $\chi_0$ son poids. Comme pour $-K_X$, on d\'ecrit les sections de $kD$ \`a l'aide de l'isomorphisme suivant:
$$\begin{array}{ccc}
\{f\in\Cbb(G/H)^{(B)}\mid \operatorname{div} (f)+kD\geq 0\} & \lra & H^0(X,kD)^{(B)} \\
  f & \lmt & fs^k.
\end{array}$$
De la proposition \ref{cardivample}, on d\'eduit alors que les poids de $H^0(X,kD)^{(B)}$ sont les caract\`eres situ\'es dans l'enveloppe convexe des $k(\chi_0 -\chi_\mathcal{C})$, o\`u les $\chi_\mathcal{C}$ sont les caract\`eres d\'efinis par le diviseur de Cartier $D$ lorsque $\mathcal{C}$ d\'ecrit l'ensemble des c\^ones colori\'es maximaux de $X$. Le polytope $Q_D^*$ est donc l'enveloppe convexe des caract\`eres $\chi_0 -\chi_\mathcal{C}$.

En particulier, les sommets de $Q_D^*$ sont entiers (c'est-\`a-dire dans $M$) et sont en bijection avec les $G$-orbites ferm\'ees de $X$, ou encore avec les plongements simples de $X$.

Par contre, si $X$ est seulement sph\'erique, il se peut que  $Q_D^*$ ait des sommets non entiers; mais les $G$-orbites ferm\'ees $Y$ de $X$ sont associ\'es \`a certains sommets entiers de  $Q_D^*$. En effet, $H^0(Y,D\cap Y)$ est un $G$-module simple, et son plus grand poids est un sommet de $Q_D^*$, car les $G$-orbites de $X$ correspondent \`a des faces de $Q_D^*$ \cite[5.3]{Br97b}.
\begin{lem}
Lorsque $X$ est horosph\'erique, avec les notations ci-dessus, les assertions suivantes sont \'equivalentes:\\
(1) $D$ est tr\`es ample,\\
(2) pour tout sommet $\chi$ de $Q_D^*$, le mono\"ide $M\cap\Rbb_+(Q_D^*-\chi)$ est engendr\'e par $M\cap (Q_D^*-\chi)$.
\end{lem}
Le mono\"ide $M\cap\Rbb_+(Q_D^*-\chi)$ est l'ensemble des \'el\'ements de $M$ qui sont dans le c\^one engendr\'e par $Q_D^*$ au point $\chi$.
\begin{center}
\includegraphics[width=13cm]{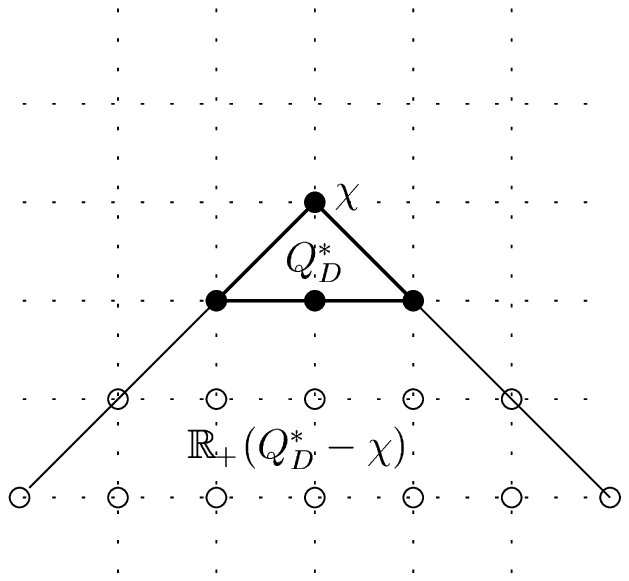}
\end{center}
\begin{rem}
Dans \cite{EW91}, G. Ewald et U. Wessels montrent que la condition (2) est v\'erifi\'ee pour tout polytope de la forme $(n-1)Q$, o\`u $Q$ est un polytope \`a sommets entiers de dimension~$n$. Or on a $Q^*_{kD}=kQ_D^*$ pour tout entier $k>0$, donc la premi\`ere partie du th\'eor\`eme \ref{tresample} se d\'eduit du lemme ci-dessus.

De plus, lorsque $X$ est sph\'erique, on a encore le r\'esultat suivant (dont la d\'emons\-tration est identique \`a celle du lemme): si la condition (2) est v\'erifi\'ee pour tout sommet entier de $Q_D^*$ alors $D$ est tr\`es ample. Ainsi la premi\`ere partie du th\'eor\`eme \ref{tresample} reste vraie dans le cas sph\'erique.
\end{rem}

\begin{proof}[D\'emonstration du lemme]
Sur une vari\'et\'e sph\'erique, tout diviseur ample est engendr\'e par ses sections globales (cons\'equence imm\'ediate de \cite[th.3.3]{Br89}), donc $D$ d\'efinit un morphisme $G$-\'equivariant
$$\begin{array}{cccc}
\phi_D : & X & \lra & \Pbb(H^0(X,D)^*) \\
  & x & \lmt & [s \mapsto s(x)].
\end{array}$$
De plus, $D$ est tr\`es ample si et seulement si $\phi_D$ est une immersion ferm\'ee, ou encore, si et seulement s'il existe un recouvrement fini de $X$ par des ouverts affines $X_{s_i}$ tels que l'application $\Cbb[\Pbb(H^0(X,D)^*)_{s_i\neq 0}]\lra \Cbb[X_{s_i}]$ soit surjective.

Les plongements simples $X_Y$ de $X$ (lorsque $Y$ parcourt les $G$-orbites ferm\'ees de $X$) recouvrent $X$, donc il suffit d'\'etudier $\phi_D$ sur chaque plongement simple $X_Y$.

Le morphisme de restriction $H^0(X,D)\lra H^0(Y,Y\cap D)$ est surjectif, car $Y$ est homog\`ene et il existe une section non identiquement nulle sur $Y$. De plus, le $G$-module $H^0(X,D)$ n'a que des multiplicit\'es $1$ dans sa d\'ecomposition en $G$-modules simples. On a ainsi une projection canonique de $H^0(X,D)^*$ dans $H^0(Y,D\cap Y)^*$  qui d\'efinit alors une application rationnelle $\pi$ de $\Pbb(H^0(X,D)^*)$ dans $\Pbb(H^0(Y,D\cap Y)^*)$. En fait $\phi_D(X_Y)$ est inclus dans le domaine de d\'efinition de $\pi$, car $X_Y=\{x\in X\mid\overline{G.x}\supset Y\}$. Notons $\psi$ la compos\'ee de $\phi_D$ et $\pi$. On a alors le diagramme suivant:
$$\xymatrix{
    X_Y\, \ar[rr]^-{\phi_D} \ar[rrd]^-{\psi}  && \Pbb(H^0(X,D)^*)\, \ar@{-->}[d]^\pi\\
     &&\Pbb(H^0(Y,D\cap Y)^*)
  }$$

On a vu au d\'ebut de cette partie que les plongements simples de $X$ sont en bijection avec les sommets de $Q_D^*$ . Soit $\chi_Y$ le sommet de $Q_D^*$ correspondant au plongement simple $X_Y$ de c\^one colori\'e $(\mathcal{C},\mathcal{F})$. Si $$W:=\Pbb(H^0(Y,D\cap Y)^*)_{v_{\chi_Y}\neq 0},$$ alors $\psi^{-1}(W)$ est un ouvert affine $B$-stable de $X_Y$; c'est donc
l'unique ouvert affine $B$-stable $X_0$ de $X_Y$ \cite[ch.3]{Kn91}.

De plus  $X_Y$ est recouvert par un nombre fini de translat\'es $g.X_0$, $g\in G$, et $$\Cbb[X_{0}]^{(B)}=\{f\in\Cbb(G/H)^{(B)} \mbox{ de poids dans } \check{\mathcal{C}}\cap M\}$$ o\`u $\check{\mathcal{C}}$ est le c\^one dual de $\mathcal{C}$ dans $M$.

On a aussi $$X_0\simeq R_u(P)\times Z$$ o\`u $P$ est un sous-groupe parabolique contenant $B$ et $Z$ est une vari\'et\'e sph\'erique affine sous l'action du sous-groupe de Levi $L$ de $P$ contenant $T$ \cite[chap. 1.4]{Br97b}. De plus \cite[chap. 3.3]{Br97b}, $$\Cbb[Z]^{(B\cap L)}=\{f\in\Cbb(Z)^{(B\cap L)} \mbox{ de poids dans } \check{\mathcal{C}}\cap M\}.$$

Remarquons maintenant que $\check{\mathcal{C}}\cap M$ n'est autre que le mono\"ide $M\cap\Rbb_+(Q_D^*-\chi_Y)$. On en d\'eduit donc que $\phi_D$ est une immersion ferm\'ee si et seulement si ce mono\"ide est engendr\'e par $M\cap (Q_D^*-\chi_Y)$, et cela pour tout $Y$.
\end{proof}

Pour d\'emontrer la deuxi\`eme partie du th\'eor\`eme \ref{tresample}, reprenons les notations de la d\'emonstra\-tion du lemme pr\'ec\'edent. Lorsque $X$ est une vari\'et\'e horosph\'erique localement factorielle, le c\^one $\mathcal{C}$ est engendr\'e par une base de $N$ d'apr\`es la proposition \ref{locfac}. Donc $\check{\mathcal{C}}$ est engendr\'e par une base de $M$. Et comme $Q_D^*$ est \`a sommets entiers, $Q_D^*-\chi_Y$ contient alors cette derni\`ere base. Ainsi la condition (2) du lemme est v\'erifi\'ee, et $D$ est tr\`es ample.

\end{document}